\documentclass[10pt]{article}

\usepackage{graphicx}
\usepackage{amscd}
\usepackage{amsmath}
\usepackage{amsfonts}
\usepackage{amssymb}
\usepackage{array}
\usepackage{theorem}

\newtheorem{theorem}{Theorem}[section]
\newtheorem{lemma}[theorem]{Lemma}
\newtheorem{corollary}[theorem]{Corollary}
\newtheorem{proposition}[theorem]{Proposition}

{\theorembodyfont{\rmfamily} \newtheorem{definition}[theorem]{Definition}
\newtheorem{example}[theorem]{Example}
\newtheorem{remark}[theorem]{Remark}
}

\numberwithin{equation}{section}

\newcommand{\zero}{\mathbf{0}}
\newcommand{\complex}{\mathbb{C}}
\newcommand{\integers}{\mathbb{Z}}
\newcommand{\supp}{\ell}
\newcommand{\dep}{{\mathop{\mathrm{dep}}\nolimits}}
\newcommand{\End}{{\mathop{\mathrm{End}}\nolimits}}

\newcommand{\pp}{\pi^{\prime}}
\newcommand{\e}{e^{(n)}}
\newcommand{\an}{\alpha^{(n)}}
\newcommand{\can}{\Check{\alpha}^{(n)}}
\newcommand{\rn}{^{(n)}}
\newcommand{\RM}{^{(m)}}
\newcommand{\LN}{\lambda^{(n)}}
\newcommand{\mn}{\mu^{(n)}}
\newcommand{\nn}{\nu^{(n)}}
\newcommand{\pln}{\pi_{\lambda^{(n)}}}
\newcommand{\plM}{\pi_{\lambda^{(m)}}}
\newcommand{\plmnn}{\mathcal{P}(\lambda,\mu,\nu,n)}
\newcommand{\pplmnn}{\mathcal{P}^+(\lambda,\mu,\nu,n)}
\newcommand{\clmn}{c_{\lambda\mu}^{\, \nu}}
\newcommand{\sn}{\Sigma^{(n)}}
\newcommand{\sm}{\Sigma^{(m)}}
\newcommand{\pmn}{\phi_{mn}}
\newcommand{\pnm}{\phi_{nm}}
\newcommand{\cp}{c_{\lambda\mu\nu}^{\; \pi}}
\newcommand{\ou}{\overline{U}\rn}
\newcommand{\ov}{\overline{V}\rn}
\newcommand{\il}{I(\lambda_1,\lambda_2)}
\newcommand{\INL}{I\rn(\lambda_1,\lambda_2)}

\newcommand{\clarge}{c_{\lambda_1 \lambda_2 \ldots \lambda_k}^{\,\nu}}
\newcommand{\csmall}{c_{\lambda_1 \lambda_2 \lambda_3}^{\,\nu}}
\newcommand{\srr}[1][X]{\Lambda^{#1}}

\begin{document}
\title{Tensor product stabilization in Kac-Moody algebras}
\author{Michael Kleber
 and Sankaran Viswanath\thanks{Research supported by a Graduate research assistantship under NSF grant DMS-9970611}}
\date{}
\maketitle
\begin{abstract}
We consider a large class of series of symmetrizable Kac-Moody
  algebras (generically denoted $X_n$). This includes the classical
  series $A_n$ as well as others like $E_n$  whose members are of
  Indefinite type. The focus is to analyze the 
behavior of representations in the limit
  $n \rightarrow \infty$. Motivated by the classical theory of 
$A_n = sl_{n+1}\complex$, we consider tensor product decompositions
  of irreducible highest weight representations of $X_n$ and study 
how these vary
  with $n$. The notion of ``double headed'' dominant weights is
  introduced. For such weights, we show that tensor product
  decompositions in $X_n$ do stabilize, generalizing the classical  results
  for $A_n$. The main tool used is Littelmann's celebrated path model.
One can also use the stable multiplicities as structure
  constants to define a multiplication operation on a suitable
  space. We define this so called {\em stable representation ring} and
  show that the multiplication operation is associative.
\end{abstract}

\section{Introduction}

In this article, we consider series of symmetrizable Kac-Moody 
algebras (generically denoted $X_n$). Our main objective is to prove
that  decompositions of tensor products of irreducible
representations of $X_n$ ``stabilize,'' i.e, given an 
irreducible representation, its multiplicity in the tensor product
decomposition becomes constant for sufficiently large $n$. 
To construct the $X_n$, let ($X, \xi$) be a marked Dynkin diagram 
with $d$  nodes and a special node $\xi$. Assume that the generalized
Cartan matrix  of $X$ is symmetrizable. We extend $X$ by
``attaching'' the Dynkin diagram
$A_{n-d}$ (a linear string of $n-d$ nodes) to $\xi$. We denote this 
new diagram $X_n$. 
\newpage
\setlength{\unitlength}{15pt}
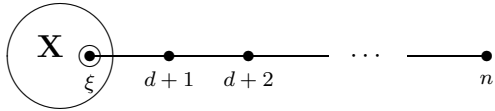
$$\begin{picture}(4,2)(1,-2)
\put(-1,0.25){\makebox(0,0){\bf \large X}}
\put(-.75,0){\circle{3}}
\put(0,0){\circle*{.25}}
\put(0,0){\circle{.5}}
\put(2,0){\circle*{.25}}
\put(4,0){\circle*{.25}}
\put(10,0){\circle*{.25}}
\put(0,0){\line(1,0){6}}
\put(10,0){\line(-1,0){2}}
\put(6.5,0){\makebox(1,0){$\cdots$}}
\put(0,-.7){\makebox(0,0){\footnotesize $\xi$}}
\put(2,-.6){\makebox(0,0){\footnotesize $d+1$}}
\put(4,-.6){\makebox(0,0){\footnotesize  $d+2$}}
\put(10,-.6){\makebox(0,0){\footnotesize  $n$}}
\end{picture}$$

The four series of finite dimensional simple Lie algebras 
$A_n, B_n, C_n, D_n$ are 
all of this form for suitable choices of $(X,\xi)$. One can 
parametrize dominant integral weights of $X_n$ by ordered pairs of
partitions. The dominant weights thus obtained are ``supported'' on both
ends of the Dynkin diagram of $X_n$. Such ``double headed'' weights
have 
been previously considered in the literature 
 \cite{bryl,hanlon,stanley,stem,benkart} in the 
context of $A_n$. Let $\mathcal{H}_2^+$ denote the set of ordered
pairs of partitions (this definition will be slightly modified in the body of this paper). 
For $\lambda, \mu \in \mathcal{H}_2^+$ we consider the corresponding 
 integrable highest weight (irreducible) 
representations $L(\LN)$ and $L(\mn)$ of $X_n$ and decompose their tensor 
product into irreducible components. 
$$ L(\LN) \otimes L(\mn) = \bigoplus \clmn (n) \, L(\nu^{(n)})$$
Here  $\clmn(n)$ denotes the multiplicity of the irreducible
representation $L(\nu^{(n)})$  in the  tensor product. 
For each fixed $\nu \in \mathcal{H}_2^+$ we prove that 
$\clmn(n)= \clmn(m)$ for all $n,m$
sufficiently 
large. We refer to this as tensor product stabilization. The main
tool used is Littelmann's path model \cite{L2} for highest weight
integrable representations of symmetrizable Kac-Moody algebras.

This result 
generalizes earlier work of R. Brylinski  \cite{bryl} on representations with double
headed highest weights for the $A_n$ case. The set of all partitions 
($\mathcal{H}_1^+$) 
can be identified with the subset of $\mathcal{H}_2^+$ of ordered
pairs 
whose second 
component is the zero partition. Our earlier association of 
double headed weights to elements of $\mathcal{H}_2^+$, 
 when restricted to $\mathcal{H}_1^+$ gives the usual
identification of partitions with dominant weights (irreducible
representations) of $A_n$. So, as a special case of our result, 
one recovers the classical $A_n$
situation, where tensor product stabilization is already implied by the 
Littlewood-Richardson rule. 

Finally, we
use the stable multiplicity values to define a new operation: the 
``stable tensor product'' on a suitably defined $\complex$ vector
space $\srr$. 
We show that this operation is associative and captures tensor product
 decompositions in the limit $n \rightarrow \infty$. 
We call $\srr$ the stable representation ring of type $X$.
In the classical
$A_n$ case, $\srr[A]$ can be viewed as the  tensor product of two
copies of the ring of symmetric functions in
infinitely many variables.

\noindent
{\bf Acknowledgements:} We would like to thank Richard Borcherds 
for encouragement and many helpful discussions. S.V would also like to thank
Peter Littelmann for his valuable input while this work was in
progress and John Stembridge for his clarifications regarding
 the type $A$ case.


\section{Formulation of the main Theorem}

\subsection{The $X_n$}\label{one}

                     We first define the series
 of symmetrizable Kac-Moody algebras that we 
will consider. Let $X$ be a Dynkin diagram in which one of the
 vertices is distinguished; we call such an object a {\em marked
 Dynkin diagram}. We assume that the associated 
generalized Cartan matrix $C(X)$ is symmetrizable; see Kac
 \cite[Chapter 4]{Kac} for background.
Let the number of nodes in $X$ be $d$. 
For convenience we number the nodes of $X$ as $1,2, \cdots, d$ 
such that the distinguished vertex is numbered
 $d$. 
For $n \geq d$, we define
$X_n$ to be the Dynkin diagram obtained from $X$ by attaching a
tail of $n-d$ nodes to the marked vertex as shown in the figure below.

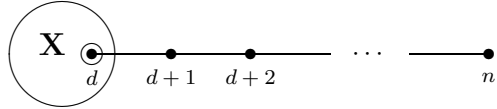
\setlength{\unitlength}{15pt}
$$\begin{picture}(4,2)(1,-1)
\put(-1,0.25){\makebox(0,0){\bf \large X}}
\put(-.75,0){\circle{3}}
\put(0,0){\circle*{.25}}
\put(0,0){\circle{.5}}
\put(2,0){\circle*{.25}}
\put(4,0){\circle*{.25}}
\put(10,0){\circle*{.25}}
\put(0,0){\line(1,0){6}}
\put(10,0){\line(-1,0){2}}
\put(6.5,0){\makebox(1,0){$\cdots$}}
\put(0,-.6){\makebox(0,0){\footnotesize $d$}}
\put(2,-.6){\makebox(0,0){\footnotesize $d+1$}}
\put(4,-.6){\makebox(0,0){\footnotesize  $d+2$}}
\put(10,-.6){\makebox(0,0){\footnotesize  $n$}}
\end{picture}$$

We ``extend'' the numbering of the nodes of $X$ to a numbering
 of the nodes of $X_n$ as in figure. Let 
$\mathfrak{g}(X_n)$ be the Kac-Moody algebra (over $\complex$) with Dynkin 
diagram $X_n$. It is  clear that $\mathfrak{g}(X_n)$ is  symmetrizable,
with generalized Cartan matrix $C(X_n)$ given by:
\begin{equation}\label{gcm}
C(X_n) = \left[
\begin{array}{c|cccc}
C(X)&&&&\\
      &-1&&&\\
\hline
    -1&2&-1&&\\
      &-1&2&\ddots& \\
      &&\ddots&\ddots&-1\\
      &&&-1&2 \\
\end{array}
\right]
\end{equation}

\begin{example}\label{eg1}
In the following diagrams, the marked vertex is the
one indicated by a circle.
\begin{enumerate}
\renewcommand{\theenumi}{\roman{enumi}}
\item If $X$ is the Dynkin diagram with a single vertex:
\setlength{\unitlength}{15pt}
\begin{picture}(1,1)
\put(0,0.5){\circle*{.25}}
\put(0,0.5){\circle{.5}}
\put(0,0){\makebox(0,0){\tiny 1}}
 \end{picture}
then $X_n$ becomes
$$\begin{picture}(4,1)
\put(0,0){\circle*{.25}}
\put(1,0){\circle*{.25}}
\put(2,0){\circle*{.25}}
\put(5,0){\circle*{.25}}
\put(0,0){\line(1,0){3}}
\put(5,0){\line(-1,0){1}}
\put(3.05,0){\makebox(1,0){...}}
\put(0,-.5){\makebox(0,0){\tiny 1}}
\put(1,-.5){\makebox(0,0){\tiny 2}}
\put(2,-.5){\makebox(0,0){\tiny 3}}
\put(5,-.5){\makebox(0,0){\tiny $n$}}
\end{picture}$$
the Dynkin diagram $A_n$.  The corresponding Lie algebra
$\mathfrak{g}(X_n) \approx sl_{n+1}(\complex)$.
We shall henceforth refer to this example as ``Type A''

\item 
\setlength{\unitlength}{15pt}
Let $X$ be the Dynkin diagram $E_6$:
$$
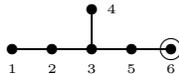
\begin{picture}(4,3)(-1,-1)
\put(0,0){\circle*{.25}}
\put(1,0){\circle*{.25}}
\put(2,0){\circle*{.25}} 
\put(2,1){\circle*{.25}}
\put(3,0){\circle*{.25}} 
\put(4,0){\circle*{.25}} 
\put(4,0){\circle{.5}}  
  \put(0,0){\line(1,0){4}}
  \put(2,0){\line(0,1){1}}
  \put(0,-.5){\makebox(0,0){\tiny 1}}
  \put(1,-.5){\makebox(0,0){\tiny 2}}
  \put(2,-.5){\makebox(0,0){\tiny 3}}
  \put(2.5,1){\makebox(0,0){\tiny 4}}
  \put(3,-.5){\makebox(0,0){\tiny 5}}
  \put(4,-.5){\makebox(0,0){\tiny 6}}
  \end{picture}
$$
For  $n \geq 6$, $X_n$ is
\setlength{\unitlength}{15pt}
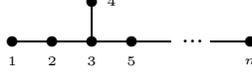
$$\begin{picture}(4,3)(-1,-1)
\put(0,0){\circle*{.25}}
\put(1,0){\circle*{.25}}
\put(2,0){\circle*{.25}} 
\put(2,1){\circle*{.25}}
\put(3,0){\circle*{.25}} 
\put(6,0){\circle*{.25}} 
  \put(0,0){\line(1,0){4}}
  \put(2,0){\line(0,1){1}}
  \put(6,0){\line(-1,0){1}}
  \put(4.05,0){\makebox(1,0){...}}
  \put(0,-.5){\makebox(0,0){\tiny 1}}
  \put(1,-.5){\makebox(0,0){\tiny 2}}
  \put(2,-.5){\makebox(0,0){\tiny 3}}
  \put(2.5,1){\makebox(0,0){\tiny 4}}
  \put(3,-.5){\makebox(0,0){\tiny 5}}
  \put(6,-.5){\makebox(0,0){\tiny $n$}}
  \end{picture}
$$
It is well known that $\mathfrak{g}(X_n)$ is a symmetrizable
Kac-Moody algebra of Finite type for $n=6,7,8$ , of Affine type
for $n=9$ and of Indefinite type for $n \geq 10$. We shall refer to this
example as ``Type E'' 

\item We can also obtain the series $B_n, C_n \;\mbox{and}\; D_n$ of
finite dimensional simple Lie algebras by choosing $X$ as follows

\begin{enumerate}
\item Type $B$: 
\begin{picture}(2,1)(-1,0)
\put(0,0.5){\circle*{.25}}
\put(1.5,0.5){\circle*{.25}}
\put(1.5,0.5){\circle{.5}}
\put(0,0.4){\line(1,0){1.5}}
\put(0,0.6){\line(1,0){1.5}}
\put(0.75,0.5){\makebox(0,0){$<$}}
\put(0,0){\makebox(0,0){\tiny 1}}
\put(1.5,0){\makebox(0,0){\tiny 2}}
\end{picture}

\item Type $C$: 
\begin{picture}(2,1)(-1,0)
\put(0,0.5){\circle*{.25}}
\put(1.5,0.5){\circle*{.25}}
\put(1.5,0.5){\circle{.5}}
\put(0,0.4){\line(1,0){1.5}}
\put(0,0.6){\line(1,0){1.5}}
\put(0.75,0.5){\makebox(0,0){$>$}}
\put(0,0){\makebox(0,0){\tiny 1}}
\put(1.5,0){\makebox(0,0){\tiny 2}}
\end{picture}

\item Type $D$:
\setlength{\unitlength}{15pt}
\begin{picture}(3,1)(-1,0)
\put(0,0.5){\circle*{.25}}
\put(1,0){\circle*{.25}}
\put(1,0){\circle{.5}}
\put(0,-0.5){\circle*{.25}}
\put(0,.5){\line(2,-1){1}}
\put(0,-.5){\line(2,1){1}}
\put(-.5,.5){\makebox(0,0){\tiny 1}}
\put(-.5,-.5){\makebox(0,0){\tiny 2}}
\put(1.5,0){\makebox(0,0){\tiny 3}}
 \end{picture}
\end{enumerate}

\item Type $F^{(1)}$: 
\setlength{\unitlength}{15pt}
\begin{picture}(3,1)(-1,0)
\put(0,0){\circle*{.25}}
\put(1,0){\circle*{.25}}
\put(2,0){\circle*{.25}}
\put(3,0){\circle*{.25}} \put(3,0){\circle{.5}}  
  \put(0,0){\line(1,0){1}}
  \put(1,.1){\line(1,0){1}}
  \put(1,-.1){\line(1,0){1}}
  \put(2,0){\line(1,0){1}}
  \put(1.5,0){\makebox(0,0){$<$}}  
  \put(0,-.5){\makebox(0,0){\tiny 1}}
  \put(1,-.5){\makebox(0,0){\tiny 2}}
  \put(2,-.5){\makebox(0,0){\tiny 3}}
  \put(3,-.5){\makebox(0,0){\tiny 4}}
  \end{picture}

\item Type $F^{(2)}$: 
\setlength{\unitlength}{15pt}
\begin{picture}(3,1)(-1,0)
\put(0,0){\circle*{.25}}
\put(1,0){\circle*{.25}}
\put(2,0){\circle*{.25}}
\put(3,0){\circle*{.25}} \put(3,0){\circle{.5}}  
  \put(0,0){\line(1,0){1}}
  \put(1,.1){\line(1,0){1}}
  \put(1,-.1){\line(1,0){1}}
  \put(2,0){\line(1,0){1}}
  \put(1.5,0){\makebox(0,0){$>$}} 
  \put(0,-.5){\makebox(0,0){\tiny 1}}
  \put(1,-.5){\makebox(0,0){\tiny 2}}
  \put(2,-.5){\makebox(0,0){\tiny 3}}
  \put(3,-.5){\makebox(0,0){\tiny 4}}
  \end{picture}

\item Type $G^{(1)}$: 
\begin{picture}(2,1)(-1.3,0)
\put(0,0.5){\circle*{.25}}
\put(1.5,0.5){\circle*{.25}}
\put(1.5,0.5){\circle{.5}}
\put(0,0.4){\line(1,0){1.5}}
\put(0,0.5){\line(1,0){1.5}}
\put(0,0.6){\line(1,0){1.5}}
\put(0.75,0.5){\makebox(0,0){$<$}}
\put(0,0){\makebox(0,0){\tiny 1}}
\put(1.5,0){\makebox(0,0){\tiny 2}}
\end{picture}

\item Type $G^{(2)}$: 
\begin{picture}(2,1)(-1.3,0)
\put(0,0.5){\circle*{.25}}
\put(1.5,0.5){\circle*{.25}}
\put(1.5,0.5){\circle{.5}}
\put(0,0.4){\line(1,0){1.5}}
\put(0,0.5){\line(1,0){1.5}}
\put(0,0.6){\line(1,0){1.5}}
\put(0.75,0.5){\makebox(0,0){$>$}}
\put(0,0){\makebox(0,0){\tiny 1}}
\put(1.5,0){\makebox(0,0){\tiny 2}}
\end{picture}

\end{enumerate}
\end{example}

\subsection{Extensible families}
For a Dynkin diagram $Y$, let $\det(Y)$ denote the determinant of the
generalized Cartan matrix of $Y$.  We allow $Y$ to be empty, in which
case $\det(Y) = 1$.
\begin{lemma}\label{detinap}
Let $X$ be a marked Dynkin diagram. 
Then, the sequence $\{ \det(X_n): n \geq d\} $ is an arithmetic progression.
\end{lemma}

\noindent
{\bf Proof:}
Let $n \geq d+2$. We can compute 
$\det(X_n)$ from  Equation~(\ref{gcm}) by expanding along the last row
of the matrix. This gives us
$$ \det(X_n) = 2 \det(X_{n-1}) - \det(X_{n-2}) \;\;\; \hfill{\Box}$$

\begin{remark} Let $\Delta$ denote the common difference of this
arithmetic progression. The argument above also works for $n=d+1$ and 
shows that $\Delta =
\det(X) - \det(X_{d-1})$ where $X_{d-1}$ denotes the Dynkin
 diagram obtained from
$X$ by deleting the distinguished vertex and all edges incident on it.
We have, for $ n  \geq d$,
\begin{equation} \label{eq1}
\det(X_n) = \det(X) + (n-d) \Delta
\end{equation}
\end{remark}

\begin{table}
\begin{center}
\setlength{\extrarowheight}{4pt}
\begin{tabular}{ | c | c | c | c | c | c |}
\hline
Type & A & B,C,D & E & $F^{(1)}$, $F^{(2)}$ & $G^{(1)}$, $G^{(2)}$ \\ \hline
$\Delta$ & 1 & 0 & -1 & -1 & -1 \\ \hline
\end{tabular}
\caption{\protect Values of $\Delta$}\label{deltatable}
\end{center}
\end{table}

\begin{definition}
The marked Dynkin diagram  $X$  is said to be {\em extensible}
if $\Delta \neq 0$, $\det(X) \neq 0$ and 
$\Delta$ is relatively prime to $\det(X)$.
\end{definition}
This technical criterion will be an assumption for all our later results.
If $X$ is extensible then Equation~(\ref{eq1}) implies that
 $\Delta$ is relatively prime to $\det(X_n)$ for all $n \geq d$. From 
Table (\ref{deltatable}) we see that Types
$A, E, F^{(i)}, G^{(i)} \, (i=1,2)$ are extensible while
Types $B,C,D$ are not.

\begin{remark}
The condition $\det(X) \neq 0$ is not an essential part of the
definition, but will be convenient for us. By Equation~(\ref{eq1}),
$\det(X_n)$ can be zero for at most one value of $n$ provided $\Delta
\neq 0$. So if $\det(X) =0$ , then $\det(X_{d+1}) \neq 0$ and 
we can replace $X$ with $X_{d+1}$
without affecting anything in the rest of this paper.
\end{remark}

\subsection{Roundup of Notation}
Most of our notation is that of Kac's book \cite{Kac}.
Let $\mathfrak{h}(X_n)$ denote the Cartan subalgebra of $\mathfrak{g}(X_n)$
and $\mathfrak{h}^*(X_n)$ denote its dual. The simple roots of 
$\mathfrak{g}(X_n)$ are denoted $\{\alpha_i^{(n)}: i=1,\cdots,n\}$. Here
$\alpha_i^{(n)}$ corresponds to the node $i$ of $X_n$ with respect to 
the node numbering mentioned in Section \ref{one}. Let 
$\Check{\alpha}_i^{(n)} 
\in \mathfrak{h}(X_n)$ be the corresponding coroot. The $(i,j)^{th}$ element
of the generalized Cartan matrix of $X_n$ is thus given by
$\alpha_j^{(n)}(\Check{\alpha}_i^{(n)})$. The root lattice of 
$\mathfrak{g}(X_n)$ is 
$$Q(X_n):=\integers \alpha_1^{(n)} 
\oplus \cdots \oplus \integers \alpha_n^{(n)} \subset \mathfrak{h}^*(X_n)$$
The weight lattice is $P(X_n) := \{ \lambda \in \mathfrak{h}^*(X_n): 
\lambda(\Check{\alpha}_i^{(n)}) \in \integers \;
\forall i=1,\cdots,n\}$
The fundamental weights 
$\omega_i^{(n)}$, $i=1,\cdots, n$ of $\mathfrak{g}(X_n)$ are elements of 
$\mathfrak{h}^*(X_n)$ which satisfy 
$\omega_i^{(n)} (\Check{\alpha}_j^{(n)}) = \delta_{ij}$. If $\det(X_n) =0$
this does not determine the  $\omega_i^{(n)}$ uniquely. In this case, we
pick them arbitrarily such that they satisfy the above condition.
We will also find it useful to index the fundamental weights ``backwards''.
We let 
$$\overline{\omega}_i^{(n)} := \omega_{n-i+1}^{(n)} \;\; i=1, \cdots, n$$
So for instance, $\omega_d^{(n)}$ is  the fundamental weight 
corresponding to the distinguished vertex of $X$ while 
$\overline{\omega}_1^{(n)}$ corresponds to the
``end'' vertex of the tail.
The set of 
dominant weights is $P^+(X_n) := \{ \lambda \in \mathfrak{h}^*(X_n): 
\lambda(\Check{\alpha}_i^{(n)}) \in \integers^{\geq 0} 
\;\; \forall i=1, \cdots, n\}$

When $\det(X_n) \neq 0$, 
$$P(X_n)=\integers \omega_1^{(n)} 
\oplus\cdots\oplus \integers \omega_n^{(n)}$$

\subsection{Double headed weights}\label{doubleheaded}

In the representation theory of $sl_{n+1}(\complex)$ (Type $A$),
dominant weights are often parametrized by partitions or equivalently
by Young diagrams. The convention is that the coefficient of the
$i^{th}$ fundamental weight $\omega_i^{(n)}$ in a given dominant
weight is the number of columns
of height $i$ in the corresponding Young diagram. A partition $\lambda$ 
with $r$ rows can thus
be thought of as defining a dominant weight $\lambda\rn$ of $A_n$ for
{\em each} $n \geq r$. We use this as motivation to similarly
parametrize weights of $X_n$. Define:
$$ \mathcal{H}_1 = \{ (x_1,x_2,\cdots) : x_i 
\in \integers \, \forall i \; and \; x_i \neq 0
 \mbox{ for only finitely many }  i \}$$
Given $x=(x_1,x_2,\cdots) \in \mathcal{H}_1$ we 
define the {\em length} of $x$ to be: $\ell(x):= \max\{i: x_i \neq 0\}$. 
The element $x \in \mathcal{H}_1$ can be used to define
a weight of $\mathfrak{g}(X_n)$ for $n \geq \supp(x)$. We let $x$ 
label the weight $x_1\omega_1^{(n)} + x_2\omega_2^{(n)} + \cdots +
 x_m \omega_m^{(n)}$ ($m=\supp(x)$) of $\mathfrak{g}(X_n)$ for $n \geq
 \supp(x)$. We also define $$ \mathcal{H}_1^+ = \{ (x_1,x_2,\cdots) : x_i 
\in \integers^{\geq 0} \, \forall i \; and \; x_i \neq 0
 \mbox{ for only finitely many }  i \}$$
By the above prescription, elements of $\mathcal{H}_1^+$ define
{\em dominant} weights of $\mathfrak{g}(X_n)$ for $n \geq \supp(x)$.
The set $ \mathcal{H}_1^+$ is also in bijection with the set of all partitions.
One  identifies $x=(x_1,x_2,\cdots) \in  \mathcal{H}_1^+$ 
with the partition $\pi$ with parts
$(x_1 + x_2 + \cdots + x_m,  x_2 + \cdots + x_m, \cdots, x_m)$. It is
easy to see that the above prescriptions  generalize that of
 the Type A situation.

There is also 
another approach to making dominant weights of different $A_n$'s 
correspond to each other. Given an ordered pair of partitions
$(\lambda,\mu)$, the convention now \cite{bryl,benkart}
 is to let the number of columns
of height $i$ in $\lambda$ be the coefficient of $\omega_i^{(n)}$ and
the number of columns of height $i$ in $\mu$ be the coefficient of
$\overline{\omega}_i^{(n)}$. Thus $\lambda$ and $\mu$ encode
information about the coefficients at the two ends of the Dynkin
diagram of $A_n$. We term such dominant weights ``double headed''.

A straightforward generalization leads to the definitions:
$\mathcal{H}_2 =\mathcal{H}_1 \times \mathcal{H}_1$ and
$\mathcal{H}_2^+ =\mathcal{H}_1^+ \times \mathcal{H}_1^+$. 
Given $x, y \in \mathcal{H}_1$, say 
$x = (x_1,x_2,\cdots), y = (y_1,y_2,\cdots)$, let
$\lambda=(x,y) \in \mathcal{H}_2$. One can 
use $\lambda$ to define a weight of $\mathfrak{g}(X_n)$ for each
$n \geq \supp(y) + \max(d, \supp({x}))$
(recall $d$ = the number of nodes in $X$) as follows:
$$\lambda^{(n)} := \sum_{i=1}^{\supp(x)} x_i \omega_i^{(n)}
+\sum_{i=1}^{\supp({y})} y_i \overline{\omega}_i^{(n)}$$
It is clear that elements of $\mathcal{H}_2^+$ define {\em dominant} weights
of $\mathfrak{g}(X_n)$. 
\noindent
We define the {\em length}:
$\supp(\lambda,X) := 
\supp(y) + \max(d, \supp(x)) $. 

In the classical Type A case, the usefulness of identifying dominant weights of
different $A_n$'s using partitions (or $\mathcal{H}_1^+$) is
apparent when studying tensor products of representations. For
instance the Littlewood-Richardson rule states that if
$V_{\lambda\rn}$ and $V_{\mu\rn}$ are the irreducible highest weight
representations correponding to partitions $\lambda$ and $\mu$, then
for large enough $n$, the tensor product $V_{\lambda\rn} \otimes
V_{\mu\rn}$  decomposes into a direct sum $\oplus
c_{\lambda\mu}^{\,\nu} V_{\nu\rn}$ . The $c_{\lambda\mu}^{\,\nu}$ here
are the Littlewood-Richardson coefficients and are independent of $n$. So
the tensor product decomposition remains essentially the same for all
large $n$. Special cases of tensor product decompositions for
  double-headed weights in Type A have been studied in  \cite{bryl}
where again one gets such a stabilization behavior for large
$n$. Double headed type $A$ weights have also been considered by
G. Benkart et al \cite{benkart} who study dimensions of correponding 
weight spaces as a function of $n$.

Analogously, given $\lambda, \mu, \nu \in 
\mathcal{H}_2^+$, we consider the irreducible representations
$L(\lambda^{(n)})$, $L(\mu^{(n)})$ and $L(\nu^{(n)})$ 
of  $\mathfrak{g}(X_n)$ with highest
weights $\lambda^{(n)}$, $\mu^{(n)}$ and $\nu^{(n)}$ respectively.
 These are all defined
provided $n$ is larger than the lengths of each of 
$\lambda$, $\mu$ and $\nu$.
The tensor product $L(\lambda^{(n)}) \otimes L(\mu^{(n)})$ is an integrable
representation of the symmetrizable Kac Moody algebra $\mathfrak{g}(X_n)$
, in category $\mathcal{O}$. It thus decomposes into a direct sum of
irreducible highest weight representations \cite[Chapter 10]{Kac}.
We let 
$c_{\lambda\mu}^{\, \nu}(n)$ denote the multiplicity of occurrence of
the representation $L(\nu^{(n)})$ in the decomposition of the tensor product 
$L(\lambda^{(n)}) \otimes L(\mu^{(n)})$. 

Note that
$c_{\lambda\mu}^{\, \nu}(n)$ is bounded above by the dimension of the
weight space $\nu\rn$ in $L(\lambda^{(n)}) \otimes L(\mu^{(n)})$.
Since all weight spaces in this representation are finite dimensional,
$c_{\lambda\mu}^{\, \nu}(n)$ is a finite number. However, if
 $\mathfrak{g}(X_n)$ is not of finite type, then there could in general be
infinitely many $\nu$ for which $c_{\lambda\mu}^{\, \nu}(n) \neq 0$.
Our main result is the following:
\begin{theorem}\label{mainthm}
Let $X$ be an extensible marked Dynkin diagram.
 Given $\lambda, \mu, \nu \in 
\mathcal{H}_2^+$, there exists a positive integer 
$N  = N(\lambda, \mu, \nu)$ 
such that
$$c_{\lambda\mu}^{\, \nu}(n) = c_{\lambda\mu}^{\, \nu}(m) \; \forall n,m \geq N$$
\end{theorem}
We denote this constant value by $\clmn(\infty)$. In general,
$N$ will  depend on $\lambda, \mu, \nu$ and $X$.
We shall prove this theorem over the course of the next two sections.

\begin{example}\label{eg2}
We consider $E_6, E_7, E_8$ with nodes numbered as in
Example~(\ref{eg1}), (ii). One has the following tensor product
decompositions:
\begin{equation}\label{tpd}
\begin{split}
E_6: \thickspace L(\omega_6^{(6)}) \otimes L(\omega_6^{(6)}) 
&= L(2\omega_6^{(6)}) \oplus L(\omega_5^{(6)}) \oplus
L(\omega_1^{(6)})\\
E_7: \thickspace L(\omega_7^{(7)}) \otimes L(\omega_7^{(7)})
&= L(2\omega_7^{(7)}) \oplus L(\omega_6^{(7)}) \oplus
L(\omega_1^{(7)}) \oplus  L(\mathbf{0}^{(7)})\\
E_8: \thickspace L(\omega_8^{(8)}) \otimes L(\omega_8^{(8)})
&= L(2\omega_8^{(8)}) \oplus  L(\omega_7^{(8)})\oplus   L(\omega_1^{(8)}) 
\oplus   L(\mathbf{0}^{(8)}) \oplus L(\omega_8^{(8)}) 
\end{split}
\end{equation}
To re-express some of this information in terms of our notations,
define the following elements of $\mathcal{H}_1^+$: 
$\mathbf{0}:=(0,0,0,\cdots)$, $\epsilon_1:=(1,0,0,\cdots)$,
$\epsilon_2:=(0,1,0,0,\cdots)$. Let $\lambda = \mu = (\zero,\epsilon_1) \in
\mathcal{H}_2^+$. Then $\lambda\rn = \mn = \overline{\omega}_1\rn =
\omega_n\rn$. For various choices of $\nu$, the values of $\clmn(n)$
for $n=6,7,8$ can be read off from Equations~(\ref{tpd}) and are given
in Table~\ref{tptable}.
\begin{table}
\begin{center}
\setlength{\extrarowheight}{4pt}
\begin{tabular}{ | c | c | c | c | c |}
\hline
$\nu$ & $\nn$ & $\clmn(6)$ & $\clmn(7)$ & $\clmn(8)$
 \\ \hline
$(\epsilon_1,\zero)$ & $\omega_1\rn$ & 1 & 1 & 1 \\ \hline
$(\zero,2\epsilon_1)$ & $2\overline{\omega}_1\rn$ & 1 & 1 & 1 \\ \hline
$(\zero,\epsilon_2)$ & $\overline{\omega}_{2}\rn$ & 1 & 1 & 1 \\ \hline
$(\zero,\zero)$ & $\zero\rn$ & 0 & 1 & 1 \\ \hline
\end{tabular}
\caption{\protect Tensor product multiplicities in $E_n$, $n=6,7,8$}
\label{tptable}
\end{center}
\end{table}
Theorem~(\ref{mainthmrestated}) will give an explicit value of $N$ for
which $\clmn(N) = \clmn(\infty)$. Using this, it will be clear that 
$ \clmn(\infty) = 0$ for $\nu=(\zero,\zero)$ and $ \clmn(\infty) = 1$
for the other three $\nu$'s in the table. 

\end{example}

\section{The {\em Number of Boxes} condition}\label{nob}

The classical Littlewood-Richardson coefficients have the property that
$c_{\lambda\mu}^{\, \nu}=0$ unless $|\lambda|+|\mu|=|\nu|$, where
$| \cdot |$ indicates the number of boxes in a Young diagram.  In this
section we give the analogous condition for double-headed weights.

Now, suppose $X$ is an extensible marked Dynkin diagram, 
and let $\lambda, \mu, \nu$ be elements of $\mathcal{H}_2^+$
then Theorem~(\ref{mainthm}) is clearly true if  
$c_{\lambda\mu}^{\, \nu}(n) = 0 $ for all large $n$. The
interesting case is when $c_{\lambda\mu}^{\, \nu}(n) \neq 0 $ for
infinitely many values of $n$. This imposes a strong compatibility
condition on $\lambda, \mu$ and $\nu$. In Type A, this condition turns
out precisely to be the number of boxes condition mentioned in the
above paragraph.

\subsection{Structure of $P(X_n)/Q(X_n)$}
First, suppose $n$ is such that $\det(X_n) \neq 0$, then it is well known that
$P(X_n)/Q(X_n)$ is a finite abelian group of order $|\det(X_n)|$.
For any $\eta \in P(X_n)$, we let $[\eta]$ denote its image in
$P(X_n)/Q(X_n)$. The following lemma motivated the extensibility criterion.
\begin{lemma}\label{cyclic}
Let $X$ be an extensible marked Dynkin diagram with $d$ nodes
and take any $n \geq d$ such that
$\det(X_n) \neq 0$. Then $P(X_n)/Q(X_n)$ is a cyclic group
with generator $[\overline{\omega}_1^{(n)}]$.
\end{lemma}
\noindent
{\bf Proof:} Since $\det(X_n) \neq 0$, $\mathfrak{h}^* (X_n)$ is spanned
over $\complex$ by the simple roots of $\mathfrak{g} (X_n)$. Consequently
$$\overline{\omega}_1^{(n)} = \sum_{i=1}^{n} k_i \alpha_i^{(n)}$$
The $k_i$'s can be  determined as follows: the entries along the 
$j^{th}$ column of $C(X_n)$ are the coefficients that one
 gets when expressing the $j^{th}$ simple root of $X_n$ in terms of the
fundamental weights. To express the $n^{th}$ fundamental
weight in terms of the simple roots, we take the inverse of $C(X_n)$ - 
the $k_i$'s are then just the entries 
along its $n^{th}$ column.
In particular 
\begin{align}
k_n &= \frac{cofactor \; of \; the \; (n,n)^{th} \; element \; of \;C(X_n)}
{\det(X_n)}\notag\\
    &= \frac{\det(X_{n-1})}{\det(X_n)}\notag
\end{align}
The extensibility of $X$ implies that
 $\det(X_n)$ and  $\det(X_{n-1})$ are relatively prime. Hence, the
smallest positive integer $c$ such that $ck_n \in \integers$ is 
$ c = |\det(X_n)|$. Thus, the order of the element 
$[\overline{\omega}_1^{(n)}]$ in $P(X_n)/Q(X_n)$
is at least $|\det(X_n)|$.  Since $P(X_n)/Q(X_n)$ has exactly 
$|\det(X_n)|$ elements, it has to be a cyclic group generated by
$[\overline{\omega}_1^{(n)}]$.

\begin{remark}
 This lemma may be false if $X$ is not extensible. For example if:
\begin{enumerate}
\item X is of Type D. Here $\Delta =0$. The group $P(D_n)/Q(D_n)$ is
  of order 4 while its subgroup generated by 
  $[\overline{\omega}_1^{(n)}]$ is only of order 2. In fact 
$P(D_n)/Q(D_n)$ fails to be a cyclic group when $n$ is even.
\item Take $X$ to be 
$$\setlength{\unitlength}{15pt}
\begin{picture}(2,1)(5,-1.3)
\put(0,0){\circle*{.30}}
\put(1.5,0){\circle*{.30}}
\put(3,0){\circle*{.30}}
\put(3,0){\circle{.50}}
\put(1.5,0){\line(1,0){1.5}}
\put(0.36,0.1){\line(1,0){0.78}}
\put(0.36,-0.1){\line(1,0){0.78}}
\put(0.40,0){\makebox(0,0){$<$}}
\put(1.10,0){\makebox(0,0){$>$}}
\put(0,-.5){\makebox(0,0){\tiny 1}}
\put(1.5,-.5){\makebox(0,0){\tiny 2}}
\put(3,-.5){\makebox(0,0){\tiny 3}}
 \end{picture}$$

\vspace{-0.3in}
This is the Dynkin diagram of affine $A_1$, extended by one more
vertex. The corresponding generalized Cartan matrix is
\begin{equation*}
C(X) = \left[
\begin{array}{rrr}
2 & -2 &0\\
-2 & 2 &-1\\
0&-1&2
\end{array}\right]
\end{equation*}
Here $\det(X) = \Delta = -2$ and hence they are not relatively
prime. In this case, the group  $P(X_n)/Q(X_n)$ has $2(n-2)$ elements
while the subgroup generated by 
  $[\overline{\omega}_1^{(n)}]$ has order $n-2$. Further
$P(X_n)/Q(X_n)$ fails to be cyclic when $n$ is even.
\end{enumerate}
\end{remark}

\noindent
The next important proposition tells us more about the images of the
fundamental weights in the groups $P(X_n)/Q(X_n)$.
\begin{proposition}\label{ai}
Let $X$ be an extensible marked Dynkin diagram with $d$ nodes and let
$\Delta$ be the common difference of $\{\det(X_n)\}_{n \geq d}$. Then,
there exists a sequence of integers $(a_i)_{i \geq 1}$ (depending
only on $X$ and the node numbering chosen) such that in $P(X_n)$
\begin{equation}\label{fundamental}
(-\Delta) \, \omega_i^{(n)} \equiv a_i \,\overline{\omega}_1^{(n)} 
\pmod{Q(X_n)}
\end{equation}
for all $i=1 ,\cdots, n$ and for all $n$ such that $\det(X_n) \neq 0$.
Further, the $a_i$'s are unique integers with this property.
\end{proposition}
\noindent
\begin{example}\label{typea}
Let $X$ be of type A: Here $\Delta = 1$ and it can be easily checked that
$a_i =i \; \forall i \geq 1$.
 We  label the vertex $i$ of 
the Dynkin diagram with the integer $a_i$ as follows:

\begin{center}
\setlength{\unitlength}{15pt}
\begin{picture}(4,1)( 1,-1)
\put(0,0){\circle*{.25}}
\put(1,0){\circle*{.25}}
\put(2,0){\circle*{.25}}
\put(5,0){\circle*{.25}}
\put(0,0){\line(1,0){3}}
\put(5,0){\line(-1,0){1}}
\put(3.05,0){\makebox(1,0){...}}
\put(5.05,0){\makebox(1,0){...}}
\put(0,-.5){\makebox(0,0){1}}
\put(1,-.5){\makebox(0,0){2}}
\put(2,-.5){\makebox(0,0){3}}
\put(5,-.5){\makebox(0,0){$n$}}
\end{picture}
\end{center}
\end{example}
Recall from Section~(\ref{doubleheaded}) that $\omega_i\rn$ is
represented by a Young diagram which is a single column of height
$i$. Thus $a_i$ ``measures'' the number of boxes in the Young diagram
corresponding to $\omega_i\rn$.

\subsection{Proof of Proposition~\eqref{ai}}
To prove Proposition~(\ref{ai}) in general, observe by 
Lemma~(\ref{cyclic}) that
for a fixed $n \geq d$ such that $\det(X_n) \neq 0$ we can find
integers
$a_1, \cdots, a_n$ such that Equation~(\ref{fundamental}) holds for
$i=1,\cdots,n$. Each of these integers is determined up to a multiple
of $\det(X_n)$. The trick is to find a single sequence $(a_i)_{i\geq
  1}$ that makes Equation~(\ref{fundamental}) hold for all $n$.

First, fix $n \geq d$ such that $\det(X_n) \neq 0$.  
Since $P(X_n)/Q(X_n)$ is cyclic with generator
$[\overline{\omega}_1^{(n)}]$, there exist $b_1, \cdots, b_n \in
\integers/(\det(X_n)) \integers$ such that
$(-\Delta)\,\omega_i^{(n)} 
\equiv b_i \, \overline{\omega}_1^{(n)} 
 \pmod{Q(X_n)}$ for $i=1,\cdots,n$. Set $R = \integers/(\det(X_n))
\integers$.
Let $b = (b_1 \; b_2 \; \cdots \; b_n)^T \in R^n$. We first obtain a
simple characterization of the $b_i$.
\begin{lemma}\label{bcharzn}
\begin{enumerate}
\renewcommand{\theenumi}{\roman{enumi}}
\item $b \in R^n$ is a solution to  $A^T b =0 \in R^n$
 where $A = C(X_n)$. Here
  we identify the elements of $A$ with their images in $R$ and treat $A$ as
an $n\times n$ matrix with entries in $R$.
\item If $x=(x_1\; x_2\; \cdots \;x_n)^T \in R^n$ is another solution to 
$A^T x =0 $, then $x$ is a multiple of $b$ .
\item $b$ is the unique element of $R^n$ such that $A^T b= 0 \in R^n$
 and $b_n = -\Delta + (\det(X_n)) \integers \in R$.
\end{enumerate}
\end{lemma}
\noindent
{\bf Proof:} 
\begin{enumerate}
\renewcommand{\theenumi}{\roman{enumi}}
\item To prove that the $i^{th}$ entry of $A^T b$ is 0 in $R$,
it is enough to show that ($i^{th}$ entry of $A^T b$)
$\overline{\omega}_1^{(n)} \equiv 0 \pmod{Q(X_n)}$. This
is because $P(X_n)/Q(X_n)$ is cyclic of order $|det(X_n)|$ with
generator $[\overline{\omega}_1^{(n)}]$. We compute:
\begin{align}
(i^{th} \text{ entry of } A^T b) \, \overline{\omega}_1^{(n)} &= 
(\sum_{j=1}^n (A^T)_{ij} b_j)\,\overline{\omega}_1^{(n)} \notag\\ 
&= \sum_{j=1}^n \alpha_i^{(n)} (\Check{\alpha}_j^{(n)})\, b_j
 \,\overline{\omega}_1^{(n)} \notag\\
&\equiv
(-\Delta)\,\sum_{j=1}^n \alpha_i^{(n)} (\Check{\alpha}_j^{(n)}) \,  
\omega_j^{(n)} \pmod{Q(X_n)}\notag
\end{align}
The last congruence just follows from the definition of the $b_j$.
We observe now that the final expression is precisely 
$(-\Delta) \alpha_i^{(n)}$.
This can be seen  by expressing $\alpha_i^{(n)}$ as a linear combination of 
the $\omega_j^{(n)}$'s and using the ``duality'' relation
 $\omega_j^{(n)}(\Check{\alpha}_k^{(n)}) = \delta_{jk}$. Clearly
$\alpha_i^{(n)} \equiv 0 \pmod{Q(X_n)}$ $\hfill{\Box}$
\item To show that any two solutions are multiples of each other, we
  will show that $A$ has an $(n-1)\times(n-1)$ minor which is a unit
  in the ring $R$. More precisely, let 
$B$ denote the principal submatrix of $A$ comprising of the
  first $n-1$ rows and columns of $A$. Observe that
$\det(B) = \det(X_{n-1})$ which is relatively prime to 
$\det(X_n)$ by the extensibility of $X$. Hence $\det(B)$ is a unit
in $\integers/(\det(X_n)) \integers$. Now
\begin{equation*}
A^T = \begin{pmatrix}
B^T & v \\
w^T & 2 
\end{pmatrix}
\end{equation*}
where $v,w \in R^{n-1}$. Since $\det(B) = \det(B^T)$ is a unit in $R$,
$(B^T)^{-1}$ exists with all its entries in
$R$. Let $C$ denote the $n\times n$ matrix
$C = \left( \begin{smallmatrix}
(B^T)^{-1} & 0 \\
0 & 1
\end{smallmatrix}
\right)$.
Then $CA^T = \left(\begin{smallmatrix}
\mathbf{I} & p \\
q^T & 2
\end{smallmatrix}
\right)$ where $p,q \in R^{n-1}$ and $\mathbf{I}$ denotes the identity
matrix of size $n-1$. We let 
$p = (p_1 \; p_2 \; \cdots \; p_n)^T$.
If $x \in R^n$ such that $A^T x = 0 \in R^n$, then 
 $C A^T x = 0$. This implies that $x_i + p_i x_n =0$ for $1 \leq
i \leq n-1$. For $x = b$ , this gives $b_i  = - p_i b_n = p_i \Delta$ since 
from its definition 
$b_n = -\Delta$. Here again, we identify all elements of $\integers$ with
  their images in $R$. Since $\Delta$ is a unit in $R$,
 $p_i = \Delta^{-1} b_i$. Substituting back , we get
\begin{equation}\label{above}
x_i = (-x_n \Delta^{-1})\, b_i \;\; \forall i \thickspace\thickspace
 \hfill{\Box}
\end{equation}
\item Follows from (i) and (ii). $\hfill{\Box}$
\end{enumerate}

We will now explicitly define the $a_i$'s. Armed with the simple
characterization of the $b_i$'s above, we will show that these $a_i$'s
satisfy Equation~(\ref{fundamental}). To construct the $a_i$'s, we
recall the notion of the dual $\Check{Y}$ of a Dynkin diagram
$Y$. This is the Dynkin diagram which corresponds to the transpose of
the generalized Cartan matrix of $Y$ i.e, $C(\Check{Y}):= C(Y)^T$.
Let us now consider $\Check{X}$ where $X$ is our given
 Dynkin diagram. For $n \geq d$ we can form $\Check{X}_n$ as before by
stipulating that the distinguished node of $\Check{X}$ be the same as that
of $X$. Clearly $\Check{X}_n$ is the dual of the Dynkin diagram $X_n$.

The Cartan subalgebra $\mathfrak{h}(\Check{X}_n)$ can be identified with
$\mathfrak{h}^*(X_n)$. The simple roots of $\Check{X}_n$ are just the simple
coroots $\Check{\alpha}_i^{(n)}$ of $X_n$ and the simple
coroots of $\Check{X}_n$ are $\alpha_i^{(n)}$. Let
$\Check{\omega}_i^{(n)} \in \mathfrak{h}(X_n)$ 
denote the 
fundamental weights
of $\Check{X}_n$ i.e, $\an_j(\Check{\omega}_i^{(n)}) = \delta_{ij}$. 
The extensibility of $X$ implies 
$\det(\Check{X}) = \det(X) \neq 0$ . Hence
 $\Check{\alpha}_i^{(n)} \, (1 \leq i \leq n)$ 
span $\mathfrak{h}^*(\Check{X}_n)
= \mathfrak{h}(X_n)$. The group $P(\Check{X})/Q(\Check{X})$ has order
$|\det(\Check{X})|$. So $\det(\Check{X}) \, \lambda \in Q(\Check{X})$
for all $\lambda \in P(\Check{X})$. We define the $a_i \, (1 \leq i
\leq d)$ by setting:
\begin{equation}\label{maineqn1}
\det(\Check{X})\, \Check{\omega}_d^{(d)} = \sum_{i=1}^d a_i 
\Check{\alpha}_i^{(d)}
\end{equation}
The argument of Lemma~(\ref{cyclic}) shows that $a_d =
\det(\Check{X}_{d-1}) = \det(X_{d-1})$.
 We define $a_i := \det(X_{i-1})$ for all $i>d$. Since
$\{ \det(X_i) : i \geq d\}$ forms an arithmetic progression, the
preceding definition of $a_i$ for $i>d$ 
and Equation~(\ref{maineqn1}) imply the following important relation:
\begin{equation}\label{maineqn2}
\det(\Check{X}_n)\, \Check{\omega}_n^{(n)} = \sum_{i=1}^n a_i 
\Check{\alpha}_i^{(n)} \;\; \forall n \geq d
\end{equation}
We claim that these $a_i$'s do our job i.e, if we fix $n \geq d$
such that $\det(X_n) \neq 0$, then
$$ (-\Delta)\, \omega_i\rn \equiv a_i \, \overline{\omega}_1\rn
 \pmod{Q(X_n)} \; \forall i=1,\cdots,n$$
It is now enough to show that
 the $a_i\, (1 \leq i \leq n)$ satisfy the condition of part (3) of
  Lemma~(\ref{bcharzn}). This is the content of the next
\begin{lemma}
\begin{enumerate}
\item Let $a =(a_1\; a_2\; \cdots \;a_n)^T \in R^n$ (usual identification). 
Then \\ $A^T a =0 \in R^n$.
\item $a_n \equiv -\Delta \pmod{\det(X_n)}$.
\end{enumerate}
\end{lemma}
\noindent
{\bf Proof:} (2) is obvious from the definition : $a_n :=\det(X_{n-1}) =
\det(X_n) - \Delta$. To prove (1), we calculate the $i^{th}$ entry of 
$A^T a$. This is equal to  
$\sum_{j=1}^n (A^T)_{ij} a_j = \sum_{j=1}^n a_j \alpha_i^{(n)} 
(\Check{\alpha}_j^{(n)}) = 
\alpha_i^{(n)} (\sum_{j=1}^n a_j\Check{\alpha}_j^{(n)})
=\alpha_i^{(n)}(\det(X_n) \Check{\omega}_n^{(n)})$, 
where the last equality uses
Equation~(\ref{maineqn2}). This final expression  is clearly 0 unless
$i=n$ in which case it is $\det(X_n)$. But $\det(X_n) = 0$ in $R$ and
we're done.$\hfill{\Box}$

For the uniqueness of the $a_i$'s observe
that if $a_i^{\prime}, \, i\geq 1$  is another such sequence for which 
Equation~(\ref{fundamental}) holds, then for each $i$, $a_i - a_i^{\prime}$
must be divisible by $\det(X_n)$ for all $n \geq i$ (for which
$\det(X_n) \neq 0$). Since
$X$ is extensible, $\Delta \neq 0$ and Equation~(\ref{eq1})
implies $|\det(X_n)| \rightarrow \infty$ as $n \rightarrow \infty$.
Hence $a_i = a_i^{\prime}$.

This finally proves Proposition~(\ref{ai}). We in fact get an explicit
method for computing the $a_i$ as well.

Equation~(\ref{maineqn2}) leads to the following additional
interpretation of the $a_i$, which we shall use later.
\begin{lemma}\label{coefofaln}
Let $X$ be an extensible marked Dynkin diagram with $d$ nodes
 and let $n \geq d$ such that $\det(X_n) \neq 0$.
Fix $i$, $1 \leq i \leq n$ and suppose $$ \omega_i^{(n)} =
\sum_{k=1}^{n} c_k \alpha_k^{(n)}$$
Then $a_i = \det(X_n) \,c_n$.
\end{lemma}
\noindent
{\bf Proof:} We have $c_n = \omega_i^{(n)} (\Check{\omega}_n^{(n)})$. 
Using Equation~(\ref{maineqn2}), we get 
\begin{align}
\det(X_n)\, c_n &= \omega_i^{(n)} (\det(X_n)\, \Check{\omega}_n^{(n)})
  &&\notag\\
                &= \omega_i^{(n)} (\sum_{j=1}^n a_j 
\Check{\alpha}_j^{(n)}) \thickspace=a_i && \hfill{\Box}\notag
\end{align}


\noindent
The next lemma and its corollary re-express the $a_i$ for $i>d$ in a
more convenient form.

\begin{lemma}\label{airev}
Let $X$ be any marked Dynkin diagram (not necessarily extensible) with $d$
 nodes.
 Let $n \geq d$ be such that $\det(X_n) \neq 0$. Then in $P(X_n)$, 
 $$\overline{\omega}_i^{(n)} \equiv i \,
 \overline{\omega}_1^{(n)} \pmod{Q(X_n)}$$
for $1 \leq i \leq (n-d+1)$.
\end{lemma}
\noindent 
{\bf Proof:} We only need to observe that if $1 \leq i \leq (n-d+1)$,
\begin{equation}\label{tailmod}
i \, \overline{\omega}_1^{(n)} - \overline{\omega}_i^{(n)} 
= \sum_{j=1}^{i-1} j \alpha_{n-i+1+j}^{(n)} \;\; \in Q(X_n) \;\; \hfill{\Box}
\end{equation}
\begin{corollary}
If $1 \leq i \leq (n-d+1)$, then $a_{n-i+1} \equiv -i\Delta 
\pmod{\det(X_n)}$.
\end{corollary}
\begin{remark}
The above corollary is also obvious from the definition of the
$a_i$. We have $a_{n-i+1} = \det(X_{n-i}) = \det(X_n) - i\Delta$.
\end{remark}
\begin{example}\label{aiE}
Type E. We indicate the $a_i$'s as labels on the Dynkin diagram.

\begin{center}
\setlength{\unitlength}{20pt}
\begin{picture}(4,3)(1,-1)
\put(0,0){\circle*{.25}}
\put(1,0){\circle*{.25}}
\put(2,0){\circle*{.25}} 
\put(2,1){\circle*{.25}}
\put(3,0){\circle*{.25}} 
\put(4,0){\circle*{.25}} 
\put(5,0){\circle*{.25}} 
\put(8,0){\circle*{.25}} 
  \put(0,0){\line(1,0){6}}
  \put(2,0){\line(0,1){1}}
  \put(8,0){\line(-1,0){1}}
  \put(6.05,0){\makebox(1,0){...}}
\put(8.05,0){\makebox(1,0){...}}
  \put(0,-.5){\makebox(0,0){\footnotesize 2}}
  \put(1,-.5){\makebox(0,0){\footnotesize 4}}
  \put(2,-.5){\makebox(0,0){\footnotesize 6}}
  \put(2.5,1){\makebox(0,0){\footnotesize 3}}
  \put(3,-.5){\makebox(0,0){\footnotesize 5}}
  \put(4,-.5){\makebox(0,0){\footnotesize 4}}
  \put(5,-.5){\makebox(0,0){\footnotesize 3}}
  \put(8,-.5){\makebox(0,0){\footnotesize $10-n$}}
  \end{picture}
\end{center}
\end{example}
\subsection{The $|\lambda|_X + |\mu|_X = |\nu|_X$ criterion}
\begin{definition}
If $\lambda = (x,y) \in \mathcal{H}_2$, we define our {\em number of boxes}
function $|\lambda|_X$ to be
\begin{equation}\label{defnob}
|\lambda|_X := \sum_{i=1}^{\supp(x)} a_i x_i - \Delta \sum_{i=1}^
{\supp(y)} i y_i
\end{equation}
\end{definition}
For instance, in our Type A example~\eqref{typea} above, 
$|\lambda|_{A}=  \sum_{i=1}^{\supp(x)} i x_i - \sum_{i=1}^
{\supp(y)} i y_i$. If we assume further that $y=(0,0,0,\cdots)$, then
$|\lambda|_{A} = \sum_{i=1}^{\supp(x)} i x_i$. 
If the dominant weight $\lambda^{(n)}$
(for $ n \geq \supp(\lambda)$) is represented as a Young diagram
(as in Section~(\ref{doubleheaded})), then
$|\lambda|_{A}$ is precisely the number of boxes in this Young diagram.
For general $y$, $|\lambda|_{A}$ measures the difference between the numbers
of boxes in the Young diagrams of $x$ and $y$.

Now, let $\lambda=(x,y) \in \mathcal{H}_2$ and fix $n \geq \supp(\lambda,X)$
such that $\det(X_n) \neq 0$. Consider the following element of $P(X_n)$ :
$(-\Delta)\lambda^{(n)} - |\lambda|_X \, \overline{\omega}_1^{(n)}$.
\begin{equation}\label{threeandhalf}
(-\Delta)\lambda^{(n)} - |\lambda|_X \, \overline{\omega}_1^{(n)}
= \sum_{i=1}^{\supp(x)} x_i ((-\Delta) \, \omega_i^{(n)} - 
a_i \overline{\omega}_1^{(n)}) + \sum_{i=1}^{\supp(y)} 
(-\Delta)\, y_i\, (\overline{\omega}_i^{(n)} - i \,\overline{\omega}_1^{(n)})
\end{equation}
The right hand side clearly lies in $Q(X_n)$ by Proposition~(\ref{ai})
and Lemma~(\ref{airev}). We have thus proved that
\begin{equation}\label{lambdamod}
(-\Delta) \lambda^{(n)} \equiv |\lambda|_X \, \overline{\omega}_1^{(n)} 
\pmod{Q(X_n)}
\end{equation}
Hence $|\lambda|_X$ identifies the coset of $Q(X_n)$ in $P(X_n)$ to
which $\LN$ belongs.
\begin{proposition}\label{numofboxes}
Let $X$ be extensible and $\lambda, \mu, \nu \in \mathcal{H}_2^+$. Suppose
$c_{\lambda\mu}^{\, \nu}(n) > 0 $ for infinitely many 
values of $n$ greater than
than the lengths of each of $\lambda, \mu, \nu$. Then
$$ |\lambda|_X + |\mu|_X = |\nu|_X$$
\end{proposition}
\noindent
{\bf Proof:}
Let $S=\{n:c_{\lambda\mu}^{\, \nu}(n)>0\} $. If $n \in S$, then the
representation $L(\nu^{(n)})$ of $\mathfrak{g}(X_n)$ occurs in the 
decomposition of the tensor product 
$L(\lambda^{(n)}) \otimes L(\mu^{(n)})$. In particular $\nu^{(n)}$
is a weight of this tensor product. All weights of 
$L(\lambda^{(n)}) \otimes L(\mu^{(n)})$ are congruent modulo the root lattice
$Q(X_n)$ to the weight $\lambda^{(n)} + \mu^{(n)}$. So, we must have
$\nu^{(n)} \equiv \lambda^{(n)} + \mu^{(n)} \; (\mbox{mod}\; Q(X_n))$.
Thus $(-\Delta) \,\nu^{(n)} \equiv 
(-\Delta)(\lambda^{(n)} + \mu^{(n)}) \pmod{Q(X_n)}$. 
Equation~(\ref{lambdamod}) then implies that
$$(|\lambda|_X + |\mu|_X - |\nu|_X) \, \overline{\omega}_1^{(n)} \equiv 0 
\pmod{Q(X_n)}$$
Finally, we use Lemma~(\ref{cyclic}) to conclude that
$|\det(X_n)|$ divides $|\lambda|_X + |\mu|_X - |\nu|_X$ for all $n \in S$.
Since $X$ is extensible, $\Delta \neq 0$ and
$|\det(X_n)| \rightarrow \infty$ as $n \rightarrow \infty$. 
This forces $|\lambda|_X + |\mu|_X - |\nu|_X =0 \;\; \hfill{\Box}$.

\begin{example}
We refer back to Example~(\ref{eg2}) and keep the same notation
here. From the definition, it is easy to see that for $\lambda
=\mu = (\zero,\epsilon_1)$, we have $|\lambda|_E = |\mu|_E =
1$. Similarly when  $\nu$  is one of $(\epsilon_1,\zero),
(\zero,2\epsilon_1)$ or $(\zero,\epsilon_2)$, $|\nu|_E =2 = 
|\lambda|_E + |\mu|_E$ while for $\nu = (\zero,\zero)$, $|\nu|_E =0$.
Proposition~\eqref{numofboxes} now implies that for $\nu = (\zero,\zero)$,
$\clmn(n) =0 $ eventually, as was stated before.
\end{example}

\section{Littelmann paths and the proof of the main theorem}

\subsection{The notion of {\em depth}}
Let $X$ be an extensible marked Dynkin diagram.
In light of Proposition~(\ref{numofboxes}), we now consider 
$\lambda, \mu, \nu \in \mathcal{H}_2^+$ such that 
$|\lambda|_X + |\mu|_X = |\nu|_X$. Let $\gamma = \lambda + \mu - \nu
\in \mathcal{H}_2$. 
Let $\gamma = (x,y)$ with 
$x, y \in \mathcal{H}_1$ and  $l:=\max(\supp(x),d)$ , $r :=
 \supp(y)$. Thus $\supp(\gamma,X)= l+r$. Since $|\gamma|_X  = 0$, we
 know that $\gamma^{(n)} \in Q(X_n)$ for all $n \geq l+r$ for which
$\det(X_n) \neq 0$ i.e, $\gamma^{(n)}$ is an integral linear combination
of $\alpha_i^{(n)} \; i=1,\cdots,n$. The next proposition tells us how
the coefficients of this linear combination change as $n$
increases. This proposition allows us to define the useful notion of depth. 
At the end of this subsection, we shall also restate our main theorem
giving an explicit value for $N$.

With notation as above, we have
\begin{proposition}\label{piqis}
There exist integers $p_i \; (1\leq i \leq l-1$) , $q_j \; (1\leq j \leq r-1)$
and $s$ such that for $n \geq l+r$
\begin{equation}\label{stringofs}
\gamma^{(n)} = \sum_{i=1}^{l-1} p_i \alpha_i^{(n)} +
\sum_{i=l}^{n-r+1} s \alpha_i^{(n)} + 
\sum_{i=n-r+2}^{n} q_{n-i+1} \alpha_i^{(n)}
\end{equation}
\end{proposition}

\begin{remark}
For the case $X_n = E_n$, the figure shows these coefficients 
labeling the corresponding nodes.
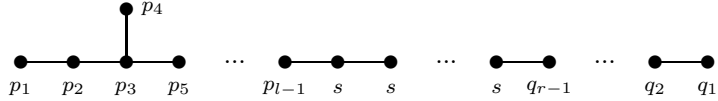
\begin{center}
\setlength{\unitlength}{20pt}
\begin{picture}(4,3)(2,-1)
\put(-1,0){\circle*{.25}}
\put(0,0){\circle*{.25}}
\put(1,0){\circle*{.25}} 
\put(1,1){\circle*{.25}}
\put(2,0){\circle*{.25}} 
\put(4,0){\circle*{.25}} 
\put(5,0){\circle*{.25}} 
\put(6,0){\circle*{.25}} 
\put(8,0){\circle*{.25}} 
\put(9,0){\circle*{.25}} 
\put(11,0){\circle*{.25}} 
\put(12,0){\circle*{.25}} 
  \put(-1,0){\line(1,0){3}}
  \put(1,0){\line(0,1){1}}
  \put(2.55,0){\makebox(1,0){...}}
  \put(4,0){\line(1,0){2}}
  \put(8,0){\line(1,0){1}}
  \put(12,0){\line(-1,0){1}}
  \put(6.55,0){\makebox(1,0){...}}
\put(9.55,0){\makebox(1,0){...}}
  \put(-1,-.5){\makebox(0,0){\footnotesize $p_1$}}
  \put(0,-.5){\makebox(0,0){\footnotesize $p_2$}}
  \put(1,-.5){\makebox(0,0){\footnotesize $p_3$}}
  \put(1.5,1){\makebox(0,0){\footnotesize $p_4$}}
  \put(2,-.5){\makebox(0,0){\footnotesize $p_5$}}
  \put(4,-.5){\makebox(0,0){\footnotesize $p_{l-1}$}}
  \put(5,-.5){\makebox(0,0){\footnotesize $s$}}
  \put(6,-.5){\makebox(0,0){\footnotesize $s$}}
  \put(8,-.5){\makebox(0,0){\footnotesize $s$}}
  \put(9,-.5){\makebox(0,0){\footnotesize $q_{r-1}$}}
  \put(11,-.5){\makebox(0,0){\footnotesize $q_2$}}
  \put(12,-.5){\makebox(0,0){\footnotesize $q_1$}}
    \end{picture}
\end{center}
\end{remark}

\noindent
Thus, as $n$ increases, the expression of $\gamma^{(n)}$
 as a linear combination of the simple roots of $X_n$ continues to
 have the same $l-1$ coefficients on the left and the same $r-1$
 coefficients on the right, while the string of $s$ 's in the middle
 grows longer.

\noindent
{\bf Proof:} We first prove the Proposition for some special choices
of $\gamma$. For $i \geq 1$, consider the following elements of
$\mathcal{H}_1$: $\sigma_i = (0,0,\cdots,-\Delta,0,0,\cdots)$ where the
$-\Delta$ occurs in the $i^{th}$ position, and $\tau_i =
(-a_i,0,0,\cdots)$. Let $\gamma_i = (\sigma_i, \tau_i) \in
\mathcal{H}_2$. Clearly $|\gamma_i|_X = 0$ for all $i$ by 
Equation~(\ref{defnob}).

Fix $i \geq 1$ and $n \geq \supp(\gamma_i, X)$ such that $\det(X_n) \neq
0$. We have  
$\gamma_i^{(n)} = (-\Delta) \omega_i^{(n)} - a_i
\overline{\omega}_1^{(n)} \in Q(X_n)$. Let 
$\gamma_i^{(n)} = \sum_{k=1}^{n} c_k \alpha_k^{(n)}$. By
Lemma~(\ref{coefofaln}) and the fact that $a_n = \det(X_{n-1})$, we
get 
\begin{equation}\label{eqone}
c_n = (-\Delta) \frac{a_i}{\det(X_n)} - 
a_i \left( \frac{\det(X_{n-1})}{\det(X_n)}\right) = -a_i
\end{equation}
For $(\max(d,i)+1) \leq j \leq n-1$,
$\gamma_i^{(n)}(\Check{\alpha}_j^{(n)}) = 0$.
But $\gamma_i^{(n)}(\Check{\alpha}_j^{(n)}) = 2 c_j - c_{j-1} -
c_{j+1}$. So
\begin{equation}\label{eqtwo}
2c_j - c_{j-1} - c_{j+1} = 0 , \; \mbox{if } \max(d,i)+1 \leq j \leq
n-1
\end{equation}
Further
\begin{equation}\label{eqthree}
\gamma_i^{(n)}(\Check{\alpha}_n^{(n)}) = -a_i = 2c_n - c_{n-1}
\end{equation}
Equations~(\ref{eqone})-(\ref{eqthree}) imply that
\begin{equation}\label{eqfour}
c_j = -a_i \text{  for } \max(d,i) \leq j \leq n
\end{equation}
We return to our general $\gamma = (x,y)$. Fix $m \geq \supp(\gamma)$
such that $\det(X_m) \neq 0$. Since $|\gamma|_X =0$, we have
\begin{equation}\label{first}
\begin{array}{ll}
\gamma^{(m)} &= (-1/\Delta) ( - \Delta \gamma^{(m)} - |\gamma|_X
\overline{\omega}_1^{(m)}) \\
&=  (-1/\Delta) \left( \sum_{i=1}^l x_i ( -\Delta \, \omega_i^{(m)} - a_i 
\overline{\omega}_1^{(m)}) - \Delta \sum_{i=1}^r y_i
(\overline{\omega}_i^{(m)} - i \, \overline{\omega}_1^{(m)}) \right) \\
&= (-1/\Delta) \sum_{i=1}^l x_i \gamma_i^{(m)} +  \sum_{i=1}^r y_i
(\overline{\omega}_i^{(m)} - i \, \overline{\omega}_1^{(m)})
\end{array}
\end{equation}

\noindent
Now if $ \sum_{i=1}^l x_i \gamma_i^{(m)} = \sum_{j=1}^n c_j
\alpha_j^{(m)}$, then Equation~(\ref{eqfour}) implies that
$c_j = - \sum_{i=1}^l x_i a_i$ for $l \leq j \leq
m$. Further, Equation~(\ref{tailmod}) implies that $\sum_{i=1}^r y_i
(\overline{\omega}_i^{(m)} - i \overline{\omega}_1^{(m)})$ is a linear
combination of $\alpha_j^{(m)}$ for $m-r+2 \leq j \leq m$. These two
observations together with Equation~(\ref{first}) mean that if
$\gamma^{(m)} = \sum_{j=1}^n d_j \alpha_j^{(m)}$, then 
$d_j = (1/\Delta)(\sum_{i=1}^l x_i a_i)$ for
$l \leq j \leq m-r+1$. We note that this implies
$(1/\Delta)(\sum_{i=1}^l x_i a_i) \in \integers$.

Define $p_i = d_i$ for $1 \leq i \leq l-1$, $q_j = d_{m-j+1}$ for
$1 \leq j \leq r-1$ and $s=(1/\Delta)(\sum_{i=1}^l x_i a_i)$. For 
$n \geq \supp(\gamma,X) = l+r$ define
$$ \mu_n = \sum_{i=1}^{l-1} p_i \alpha_i^{(n)} +
\sum_{i=l}^{n-r+1} s \alpha_i^{(n)} + 
\sum_{i=n-r+2}^{n} q_{n-i+1} \alpha_i^{(n)} \in Q(X_n)$$
By definition, $\mu_m = \gamma^{(m)}$. Now for $1 \leq i \leq l$,
$\mu_n(\Check{\alpha}_i^{(n)})$ only depends on the values $p_i$,
$p_j$ for $j$ running over all neighbors of the node $i$ in $X_n$
and possibly on $s$ (if $i=l$ or $l-1$). Thus
$\mu_n(\Check{\alpha}_i^{(n)})$
is independent of $n$. Similarly, $\mu_n(\Check{\alpha}_{n-j+1}^{(n)})$
is independent of $n$ for $1 \leq j \leq r$. Further 
$\mu_n(\Check{\alpha}_i^{(n)}) =0$ for $l+1 \leq i \leq n-r$. 
These facts combined with  $\mu_m = \gamma^{(m)}$ gives us that
 $\mu_n = \gamma^{(n)}$ for all $n \geq l+r$. $\hfill{\Box}$

\begin{definition}
If $\gamma$ is any element of $\mathcal{H}_2$ such that $|\gamma|_X
=0$, it is clear that Proposition~(\ref{piqis}) still holds. 
We shall call the number $s$ that occurs in Proposition~(\ref{piqis})
the {\em depth} of $\gamma$. We write 
$$\dep(\gamma) := s = (1/\Delta) \sum_{i=1}^l x_i a_i = \sum_{j=1}^r
j\, y_j $$
\end{definition}
The last equality follows from $|\gamma|_X = 0 $.

\begin{lemma}
Let $\lambda, \mu, \nu \in \mathcal{H}_2^+$ be such that 
$|\lambda|_X + |\mu|_X = |\nu|_X$. Suppose 
$c_{\lambda\mu}^{\,\nu}(n) > 0 $ for some $n \geq
 \supp(\lambda+\mu-\nu,X)$, then $\dep(\lambda+\mu-\nu) \geq 0$.
\end{lemma}
\noindent
{\bf Proof:} We have $\lambda^{(n)} + \mu^{(n)} - \nu^{(n)} \in
Q^+(X_n)$. So if $\lambda^{(n)} + \mu^{(n)} - \nu^{(n)} =
\sum_{i=1}^n d_i \alpha_i^{(n)}$, then all the $d_i \geq 0$. 
By Proposition~(\ref{piqis}), we now conclude that
$\dep(\lambda+\mu-\nu) \geq 0 \hfill{\Box}$

We restate our main Theorem~(\ref{mainthm}) for the case $|\lambda|_X
+ |\mu|_X = |\nu|_X$ giving an explicit value for $N$.
\begin{theorem}\label{mainthmrestated}
Let $X$ be an extensible marked Dynkin diagram
and $\lambda,\mu,\nu \in \mathcal{H}_2^+$ such
that $|\lambda|_X + |\mu|_X = |\nu|_X$. Let $\gamma = \lambda + \mu
-\nu \in  \mathcal{H}_2$ and $N = \supp(\gamma,X) + 2 \,\dep(\gamma)$. 
Then $\clmn(m) = \clmn(n)$ for all $n,m \geq N$. 
We denote this constant value by $\clmn(\infty)$ as before.
\end{theorem}
We shall prove this theorem in the next few subsections. For the rest
of this section, $\lambda, \mu, \nu, \gamma, N$ will be as in the
statement of this Theorem. By (\ref{piqis}) we know that
$$\gamma^{(n)} = \sum_{i=1}^{l-1} p_i \alpha_i^{(n)} +
\sum_{i=l}^{n-r+1} s \alpha_i^{(n)} + 
\sum_{i=n-r+2}^{n} q_{n-i+1} \alpha_i^{(n)}$$
where $s=\dep(\gamma)$. Here $l, r, p_i, q_j$ are all as in
Proposition~\eqref{piqis}. 

\subsection{The path model}
As a first step in proving Theorem~(\ref{mainthmrestated}) we will need an explicit
expression for $c_{\lambda\mu}^{\,\nu}(n)$ given by Littelmann's path
model  \cite{L2}. We recall the relevant notions
here.

Let $\Pi^{(n)}$ denote the set of all piecewise linear paths $\pi:
[0,1] \rightarrow  \mathfrak{h}^*(X_n)$ such that $\pi(0)=0$. We
identify paths that are reparametrizations of each other. For each 
simple root $\alpha_i^{(n)} \,(1 \leq i \leq n)$, we define a 
{\em lowering operator}  $f_i^{(n)}$ and a {\em raising operator}
$e_i^{(n)}$ on $\integers \Pi$, the free $\integers$ module with basis
$\Pi$. Given $\pi \in \Pi^{(n)}$, let 
$\pi_i(t) = \pi(t)(\Check{\alpha}_i^{(n)})$ for $0 \leq t \leq 1$.
 We consider the function $a: [0,1] \rightarrow [0,1]$ defined by
$a(t) = \min\{1, \pi_i(s) - m_i | t \leq s \leq 1\}$, where $m_i =
 \min\{\pi_i(t)|0 \leq t \leq 1\}$.  Note that $a$ is an increasing
 function.
If $a(1) < 1$ , $f_i^{(n)} \pi
 :=0$. Otherwise, $f_i^{(n)} \pi$ is the path defined by
\begin{equation}\label{fi}
f_i^{(n)} \pi (t) := \pi (t) - a(t) \alpha_i^{(n)}
\end{equation}
So if $f_i^{(n)} \pi \neq 0$, then 
\begin{equation}\label{filower}
f_i^{(n)} \pi (1) = \pi(1) -
\alpha_i^{(n)}
\end{equation}
 Thus $f_i\rn$ lowers the endpoint of the path $\pi$
by $\alpha_i\rn$.

Similarly we consider the increasing function 
$b: [0,1] \rightarrow [0,1]$ with
$b(t) = \max\{0, 1 - (\pi_i(s) - m_i) | 0 \leq s \leq t\}$. If $b(0) > 0$,
we set $e_i^{(n)} \pi =0$ and otherwise 
\begin{equation}\label{ei}
e_i^{(n)} \pi(t) := \pi (t) + b(t) \alpha_i^{(n)}
\end{equation}
If $e_i^{(n)} \pi \neq 0$, then $e_i^{(n)} \pi (1) = \pi(1) +
\alpha_i^{(n)}$. For a more ``geometric'' description of the action of
the lowering and raising operators, see Littelmann  \cite{L1,L2,L3}.

\begin{remark}\label{me1}
We consider the following situation which will occur often. If
$\pi_i(t)$ is itself an increasing function with $\pi_i(1)=1$, then 
from the definition, we get $a(t)=\pi_i(t)$.
\end{remark}

To obtain the value of $c_{\lambda\mu}^{\,\nu}(n)$,
 we first consider the straight line path $\pi_{\lambda^{(n)}} \in \Pi^{(n)}$
 defined by $\pln(t) = t \LN$ for $t \in [0,1]$. 
The set of all paths that can be
 obtained by repeated action of the lowering operators on $\pln$ is
 called the set of Lakshmibai-Seshadri (L-S) paths of shape $\LN$. Let
 $$\mathcal{P}(\lambda, \mu, \nu, n) := \{ 
\text{ L-S paths of shape }
 \LN \text{ whose endpoint is } \nn - \mn\}$$
If $\pi=f_{i_k}\rn \cdots
 f_{i_2}\rn f_{i_1}\rn (\pln)$ is an element of $\plmnn$, then clearly
 Equation~(\ref{filower}) implies that $\sum_{j=1}^k \alpha_{i_j}\rn =
\LN + \mn - \nn$. A path $\pi \in \plmnn$  is said to be {\bf $\mn$
 dominant} if the translated path $\mn + \pi(t)$ lies completely in
 the dominant Weyl chamber of $\mathfrak{h}^*(X_n)$. Let 
$$\pplmnn :=\{ \pi \in \plmnn \,:\, \pi \text{ is }\mn \text{ dominant}\}$$
Littelmann's tensor product decomposition
 formula \cite{L2} now  states that the number of elements in 
$\pplmnn$ is the value of $\clmn(n)$.
\begin{theorem}\label{litt}
{\rm (Littelmann)}  $\clmn(n)= \# \pplmnn$ 
\end{theorem}

This theorem will be the main tool in our proof of Theorem~(\ref{mainthmrestated}).

\subsection{}
In light of Theorem~(\ref{litt}), one needs to analyze the 
set $\pplmnn$ better. In this subsection, we introduce certain special
lowering operators. It will turn out that paths in $\pplmnn$ can be
obtained by repeated application of just these special lowering
operators on $\pi_{\lambda\rn}$. This fact will imply our main
theorem~\eqref{mainthmrestated}.

We first  consider a larger
set of paths. Let 
$$V\rn := \{i: l+s < i < (n-r+1)-s\}$$ and
 $$\overline{V}\rn := \{i: l+s \leq i \leq (n-r+1)-s\}$$
Let $\Sigma\rn \subset \Pi\rn$ be 
$$ \Sigma\rn := \{ \eta \in \Pi\rn \,|\, \eta(t)(\Check{\alpha}_i\rn) =0
\,\forall t\in[0,1]; \, \forall i \in V\rn\}$$
i.e, $ \Sigma\rn$ is the set of paths that are ``supported'' on 
$l+s$ nodes on the left and $r+s$ nodes on the right.
Now, $\integers \Sigma\rn$ will no longer be closed under the action
of all the lowering operators. We will show below that there are still 
many lowering operators and certain compositions of them that 
preserve $\integers\Sigma\rn$. 
Let 
\begin{align}
 g\rn &:= f_{n-r+1-s}\rn \cdots f_{l+s+1}\rn f_{l+s}\rn &\text{ and }&
 &h\rn &:= f_{l+s}\rn f_{l+s+1}\rn\cdots f_{n-r+1-s}\rn \notag
\end{align}
Clearly $g\rn,
h\rn \in \End \,\integers\Pi\rn$. 
\begin{lemma}\label{sigmainv}
Suppose $\eta \in \sn$ . Then
\begin{enumerate}
\item Let $j \not\in \overline{V}\rn$.
If $f_j\rn \eta \neq 0$, then
$f_j\rn \eta \in \sn$.
\item  If $g\rn \eta \neq 0$, then
$g\rn \eta \in \sn$.
\item If $h\rn \eta \neq 0$, then
$h\rn \eta \in \sn$.
\end{enumerate}
\end{lemma}
\noindent
{\bf Proof:} (1) By Equation~(\ref{fi}),  
$f_j^{(n)} \eta (t) = \eta(t) - a(t) \alpha_j^{(n)}$.
If $i \in V\rn$, then the nodes $i$ and $j$ are not connected by a line
in the Dynkin diagram $X_n$. Hence
$\alpha_j\rn(\Check{\alpha}_i\rn)=0$. This together with $\eta \in
  \sn$ implies that $f_j\rn \eta \in \sn$.\\
(2) Suppose $(f_{l+s}\rn \eta)(t) = \eta(t) - a(t) \alpha_{l+s}\rn$,
  then 
\begin{equation*}
\begin{array}{ll}
(f_{l+s}\rn \eta)(t) (\Check{\alpha}_{l+s+1}\rn) &=
  \eta(t)(\Check{\alpha}_{l+s+1}\rn) - 
a(t) \alpha_{l+s}\rn(\Check{\alpha}_{l+s+1}\rn)\\ 
&= a(t)
\end{array}
\end{equation*}
Since $g\rn \eta \neq 0, f_{l+s}\rn \eta \neq 0$. Hence 
$a(t)$ is an increasing function with $a(1)=1$. By
Remark~(\ref{me1}) we have 
$ (f_{l+s+1}\rn f_{l+s}\rn \eta)(t) = \eta(t) - a(t) \alpha_{l+s}\rn -
  a(t) \alpha_{l+s+1}\rn $.
Continuing this process, we have
\begin{equation}\label{gn}
(g\rn\eta)(t) = (f_{n-r+1-s}\rn \cdots f_{l+s+1}\rn f_{l+s}\rn \eta)(t)
=\eta(t) - a(t) \sum_{j=l+s}^{n-r+1-s} \alpha_j\rn
\end{equation}
But $\eta \in \Sigma\rn$ and 
$(\sum_{j=l+s}^{n-r+1-s} \alpha_j\rn)(\Check{\alpha}_i\rn)=0$ for all
$i \in V\rn$. Hence, $g\rn \eta \in \Sigma\rn$ too. The proof of (3) is
analogous. $\hfill{\Box}$

The definition of $\sn$ makes it clear that $\sn$ and $\Sigma^{(m)}$
are in some sense the same, since the paths in both sets are basically
just supported on $l+r+2s$ nodes. To make this more precise, we
define maps $\phi_{nm}: \sn \rightarrow \sm$ for all $n,m \geq
l+r+2s$ as follows: Take $\eta \in \sn$. Since
$\eta(t)(\Check{\alpha}_i\rn) =0$ for all $i \in V\rn$, we can write
$$\eta(t) = \sum_{i=1}^{l+s} d_i(t) \omega_i\rn + \sum_{j=1}^{r+s}
\tilde{d}_j(t) \overline{\omega}_j\rn$$
We define
$$\pnm(\eta)(t) := \sum_{i=1}^{l+s} d_i(t) \omega_i\RM + \sum_{j=1}^{r+s}
\tilde{d}_j(t) \overline{\omega}_j\RM$$
Clearly $\pnm$ and $\pmn$ are inverses of each other and set up
bijections between the sets $\sn$ and $\sm$.

The following lemma ensures that these bijections also respect the
actions of the special lowering operators introduced above. We let
$\overline{f}_i\rn$ denote the lowering operator $f_{n-i+1}\rn$.
\begin{lemma}\label{propofsn}
Let $m,n \geq l+r+2s$ and $\eta \in \sn$.
\begin{enumerate}
\item If $1 \leq i < l+s$ then $\pnm(f_i\rn \eta) = f_i\RM
  \pnm(\eta)$.
\item If $1 \leq j < r+s$ then $\pnm(\overline{f}_j\rn \eta) = 
\overline{f}_j\RM  \pnm(\eta)$.
\item $\pnm(g\rn \eta) = g\RM
  \pnm(\eta)$.
\item $\pnm(h\rn \eta) = h\RM
  \pnm(\eta)$.
\end{enumerate}
All these equalities also hold if some of the paths involved become
0. We define $\pnm(0)=0$
\end{lemma}
\noindent
{\bf Proof:} (1) and (2) follow  from the definitions of the lowering
operators and $\pnm$. 
For (3), suppose $\eta(t) = \sum_{i=1}^{l+s} d_i(t) \omega_i\rn +
\sum_{j=1}^{r+s} \tilde{d}_j(t) \overline{\omega}_j\rn$, then
 Equation~(\ref{gn}) implies
that $$(g\rn \eta)(t)= \eta(t) - a(t) \sum_{j=l+s}^{n-r+1-s}
\alpha_j\rn$$
But $\sum_{j=l+s}^{n-r+1-s}\alpha_j\rn = \omega_{l+s}\rn + \overline{\omega}_{r+s}\rn$.
Hence $$(g\rn \eta)(t)= \eta(t) - a(t) (\omega_{l+s}\rn +
\overline{\omega}_{r+s}\rn)$$ It is easy to see that if we replace $n$ by
$m$ and $\eta$ by $\pnm(\eta)$ throughout, then the above argument
still holds, showing that $(g\RM \pnm(\eta))(t)= \pnm(\eta) - a(t)
(\omega_{l+s}\RM + \overline{\omega}_{r+s}\RM)$. 
The proof of (4) is similar. $\hfill{\Box}$

Since these special lowering operators seem to be natural in our
setting, we next consider the subset of L-S paths of shape 
$\lambda\rn$ which are obtained by repeated actions of only these
special lowering operators. More precisely, define

\vspace{0.2in}
\noindent
$\mathcal{P}^0(\lambda,\mu,\nu,n):=$ set of L-S paths
$\pi$ of shape $\LN$, with $\pi(1) = \nn - \mn$ such that $\pi$ can be
obtained by the action of the operators $\{ f_i\rn:1\leq i<l+s\} \cup
\{\overline{f}_j\rn: 1\leq j<r+s\} \cup \{g\rn,h\rn\}$ on $\pln$.

\vspace{0.1in}
We then have the following: 
\begin{lemma}\label{posubpo}
\begin{enumerate}
\item $\mathcal{P}^0(\lambda,\mu,\nu,n) \subset \sn$.
\item 
$\pnm(\mathcal{P}^0(\lambda,\mu,\nu,n)) \subset 
\mathcal{P}^0(\lambda,\mu,\nu,m) $.
\end{enumerate}
\end{lemma}
\noindent
{\bf Proof:} (1) Since the path $\pln \in \sn$, acting on it by the
special lowering operators still gives us a path in $\sn$
(by Lemma~(\ref{sigmainv})).\\
(2) If the path $\pi$ is obtained
by the action of the operators  $f_i\rn, \overline{f}_j\rn,
g\rn$ and $h\rn$ on $\pln$, Lemma~(\ref{propofsn}) implies that 
$\pnm(\pi)$ is obtained by the action of the
corresponding operators  $f_i\RM, \overline{f}_j\RM,
g\RM$ and $h\RM$ on $\pnm(\pln)$. But $\pnm(\pln)=\plM$ since the
support of $\lambda$ is a subset of the first $l$ and last $r$ nodes.
Further, since the endpoint of $\pi$ is 
$\nn - \mn$, the endpoint of $\pnm(\pi)$ is clearly $\nu\RM -
\mu\RM$. Thus  $\pnm(\pi) \in \mathcal{P}^0(\lambda,\mu,\nu,m)$. 
$\hfill{\Box}$

Clearly $\mathcal{P}^0(\lambda,\mu,\nu,n)$ and
$\mathcal{P}^+(\lambda,\mu,\nu,n)$ are both subsets of 
 $\mathcal{P}(\lambda,\mu,\nu,n)$. The next important proposition 
relates these subsets.

\begin{proposition}\label{thecrux}
$\mathcal{P}^+(\lambda,\mu,\nu,n) \subset
\mathcal{P}^0(\lambda,\mu,\nu,n)$
\end{proposition}

\noindent
Before we embark upon the proof of Proposition~(\ref{thecrux}), we
state a corollary which implies our main Theorem~(\ref{mainthmrestated}).

\begin{corollary}\label{maincor}
$\pnm(\mathcal{P}^+(\lambda,\mu,\nu,n)) \subset 
\mathcal{P}^+(\lambda,\mu,\nu,m) $
\end{corollary}
\noindent
{\bf Proof:} Let $\pi \in
\mathcal{P}^+(\lambda,\mu,\nu,n)$. Proposition~(\ref{thecrux}) implies
that $\pi \in \mathcal{P}^0(\lambda,\mu,\nu,n)$. By 
Lemma~(\ref{posubpo}), $\pnm(\pi) \in
\mathcal{P}^0(\lambda,\mu,\nu,m)$ ; in particular $\pnm(\pi)$ is an
L-S path. We need to show that $\pnm(\pi)$ is $\mu\RM$ dominant. Since
$\pi \in \mathcal{P}^0(\lambda,\mu,\nu,n \,\subset \Sigma\rn$, write
$$\pi(t) = \sum_{i=1}^{l+s} d_i(t) \,\omega_i\rn + \sum_{j=1}^{r+s}
\tilde{d}_j(t) \,\overline{\omega}_j\rn$$ 
Let $\mu = (x,y) \in \mathcal{H}_2^+$.
Since $\pi$ is $\mn$ dominant, we have
$(\mn+\pi(t))(\Check{\alpha}_i\rn) \geq 0 \, \forall t$ for all $1
\leq i \leq n$. This is equivalent to the following conditions
\begin{enumerate}
\item $x_i + d_i(t) \geq 0 \, \forall t$; $1 \leq i \leq l+s$
\item $y_j + \tilde{d}_j(t) \geq 0 \, \forall t$; $1 \leq j \leq r+s$
\end{enumerate}
It is clear that these very same conditions imply the
fact that $\pnm(\pi)$ is $\mu\RM$ dominant. $\hfill{\Box}$

\noindent
{\bf \underline{Proof of Theorem~\eqref{mainthmrestated}}:}
Corollary~(\ref{maincor}) together with the fact that $\pnm$ and
$\pmn$ are inverse maps imply that the sets 
$\mathcal{P}^+(\lambda,\mu,\nu,n)$ and 
$\mathcal{P}^+(\lambda,\mu,\nu,m)$ are in bijection with each other,
for $m,n \geq l+r+2s$. 
We now appeal to Theorem~(\ref{litt}) to deduce 
 Theorem~(\ref{mainthmrestated}): $\clmn (n) = \clmn (m)$
provided $n,m \geq l+r+2s$. $\hfill{\Box}$

\subsection{Proof of Proposition~(\ref{thecrux})}
To prove proposition~\eqref{thecrux}, we shall start with a path $\pi
\in \pplmnn$ and construct a string of raising operators which maps
$\pi$ to $\pi_{\lambda\rn}$. These raising operators will be the
analogues of the special lowering operators introduced before.

We first state some properties of  raising operators that we will
need. We refer to Littelmann's paper  \cite{L2} for the proofs.

\begin{proposition}\label{eifacts}
Let $\eta$ be an element of $\Pi\rn$ and $1 \leq i, j \leq n$.
\begin{enumerate}
\item If the nodes $i$ and $j$ have no
edge between them i.e, $\an_i(\can_j)=0=\an_j(\can_i)$, then 
$\e_i \e_j \eta = \e_j \e_i \eta$.
\item $\e_i \eta = 0 \Leftrightarrow \eta(t)(\can_i) \geq 0 \; \forall
  t \in [0,1]$.
\item If $\e_i \eta \neq 0$, then $\min_t (\e_i \eta)(t)(\can_i) =
\min_t \eta(t)(\can_i)   + 1 $.
\item If $\e_i \eta \neq 0$, then $f_i\rn \e_i \eta = \eta$
\item If $\eta$ is an L-S path, then $\eta$ has the
  {\em integrality property} i.e, $\min_t \eta(t)(\can_i)$ is an integer
  for all $1 \leq i \leq n$.
\item If $\eta$ is an L-S path of shape $\LN$ then 
$\LN - \eta(1) \in Q^+(X_n)$.
\end{enumerate}
\end{proposition}

\noindent
We shall now prove Proposition~(\ref{thecrux}). Let $U_1 = \{ l < i<
n-r+1\}$ and $\overline{U}_1 = \{ l \leq i \leq n-r+1\}$. Assume 
$\pi \in \mathcal{P}^+(\lambda,\mu,\nu,n)$. By definition, this means
that $\mn(\can_i) + \pi(t)(\can_i) \geq 0$ for all $t$ and for all 
$ 1 \leq i \leq n$. In particular, since $\mn(\can_i) = 0 $ for all
$ i \in U_1$, we have 
\begin{equation}\label{prop-one}
\pi(t)(\can_i) \geq 0 \; \forall t \in [0,1], \, \forall i \in U_1
\end{equation}
Secondly, since $\pi(1) = \nn - \mn$, we have
\begin{equation}\label{prop-two}
\pi(1)(\can_i) = 0 \, \forall i \in U_1
\end{equation}
Properties~(\ref{prop-one}) and (\ref{prop-two}) will be important for
us. In fact we will only need these two properties of $\pi$ and the
fact that $\pi$ is an L-S path to show that
$\pi \in \mathcal{P}^0(\lambda,\mu,\nu,n)$.

Since $\pi$ is an L-S path, there exists a sequence of raising
operators which maps $\pi$ to $\pln$. Let $\e_{i_p}\cdots
\e_{i_2}\e_{i_1} \pi = \pln$. Clearly $\LN- \pi(1) = \sum_{k=1}^p
\an_{i_k}$. Pick $j$ minimal such that $i_j \in \overline{U}_1$. 

\noindent
{\bf \underline{Claim}:} $i_j 
\not\in U_1$. 

\noindent
{\it Proof:} Suppose $i_j \in U_1$. For $1 \leq k \leq
j-1$, $i_k \not\in \overline{U}_1$. Hence the nodes $i_j$ and $i_k$ of
$X_n$ do not have an edge between them. By (1) of
Proposition~(\ref{eifacts}), this implies that $\e_{i_j}$ commutes
with $\e_{i_k}$ for all $1 \leq k \leq j-1$. Thus 
$\e_{i_j}\cdots \e_{i_2}\e_{i_1} \pi = \e_{i_{j-1}}\cdots
\e_{i_2}\e_{i_1} \e_{i_j}\pi =  0 $; since
$\e_{i_j} \pi =0$  by
Property~(\ref{prop-one}) and Proposition~(\ref{eifacts}), (2).
This contradicts 
$\e_{i_p}\cdots \e_{i_2}\e_{i_1} \pi = \pln \neq 0$, proving our
claim $\Box$

So either $i_j = l$ or $i_j = n-r+1$. {\bf Case 1:} $i_j
=l$. Let $\pi^{\prime} = \e_{i_j}\cdots \e_{i_2}\e_{i_1} \pi $. Then
$\pi^{\prime}(1) = \pi(1) + \sum_{k<j} \an_{i_j} + \an_l$. By
Property~(\ref{prop-two}), $\pi^{\prime}(1)(\can_{l+1}) =
\an_l(\can_{l+1}) = -1$. Proposition~(\ref{eifacts}), (2) implies that
$\e_{l+1}\pi^{\prime} \neq 0$. Again, by definition
$(\e_{l+1}\pi^{\prime})(1) = \pi^{\prime}(1) + \an_{l+1}$. So
$(\e_{l+1}\pi^{\prime})(1)(\can_{l+2}) = -1$. Thus
 $\e_{l+2}\e_{l+1}\pi^{\prime} \neq 0$. We continue this way to
conclude that $\e_{n-r}\cdots \e_{l+2}\e_{l+1} \pp \neq 0$. Note
  that we cannot go all the way to $\e_{n-r+1}$ since
  $\pi(1)(\can_{n-r+1})$ may not be 0. We set $\pi_1 :=
\e_{n-r}\cdots \e_{l+2}\e_{l+1} \pp$. {\bf Case 2:} If $i_j = n-r+1$,
the same argument as in Case 1 proves that
$\e_{l+1}\cdots\e_{n-r-1}\e_{n-r} \pp \neq 0$. In this case, we set
$\pi_1 :=\e_{l+1}\cdots\e_{n-r-1}\e_{n-r} \pp$. 

By Proposition~(\ref{eifacts}), (4), we have just shown that $\pi$ can be obtained by repeated action of
operators from the set $\{f_i\rn: i \not\in V\rn\} \cup \{g\rn,h\rn\}$
on $\pi_1$. We recall that $V\rn = \{l+s < i < n-r+1-s\}$. Note
 that in either case, 
\begin{equation}\label{fin}
\pi_1(1) - \pi(1) \geq \sum_{i=l+1}^{n-r} \an_i
\end{equation}
Here we used the `usual' partial order on $P(X_n)$ defined by
$\alpha \geq \beta \Leftrightarrow \alpha - \beta \in Q^+(X_n)$.

Now, $\pi_1$ is still an L-S path of shape $\LN$. We define 
the sets $U_2 := \{ i: l+1 < i < n-r\}$ and $\overline{U}_2 := 
\{i: l+1 \leq i \leq n-r\}$, obtained by deleting one node from each
end of the string of nodes in $U_1$ and $\overline{U}_1$. One observes
that (a) $\pp(1) (\can_i) = 0 $ for all $i
\in U_2$ and (b) $(\sum_{i=l+1}^{n-r} \an_i) (\can_i) =0 \, \forall i
\in U_2$. Since $\pi(1) = \pp(1) + \sum_{i=l+1}^{n-r} \an_i$,
these give us
\begin{equation}\label{prop-three}
\pi_1(1)(\can_i) =0 \, \forall i \in U_2
\end{equation}
This is similar to Property~(\ref{prop-two}) of $\pi$. We now claim
that the analog of Property~(\ref{prop-one}) of $\pi$ also holds for
$\pi_1$. More precisely we have:
\begin{lemma}
\begin{equation}\label{prop-four}
\pi_1(t)(\can_i) \geq 0 \, \forall t \in [0,1], \,\forall i \in U_2
\end{equation}
\end{lemma}
\noindent
{\bf Proof:}
Let $i \in U_2$ i.e, $l+2 \leq i \leq (n-r-1)$. We only consider 
{\bf Case 1:} $\pi_1
= \e_{n-r}\cdots\e_{l+2}\e_{l+1} \pp$. The other case will follow by a similar
argument. Let $\eta = \e_{i-1}\cdots\e_{l+1} \pp$.
Then by succesively using the definitions of $\e_i, \e_{i+1}, \cdots,
\e_{n-r}$ we get
\begin{equation}
\begin{array}{ll}
\pi_1(t) &= (\e_{n-r}\cdots\e_{i+1}\e_i \eta)(t) \\
&= \eta(t) + b_i(t)\an_i + b_{i+1}(t) \an_{i+1} + \sum_{k=i+2}^{n-r}
b_k(t) \an_k
\end{array}
\end{equation}
where the $b_j:[0,1] \rightarrow [0,1]$ are increasing functions with
$b_j(0)=0$ and $b_j(1)=1$. Note that $(\e_i \eta)(t) = \eta(t) + b_i(t)
\an_i$. Similarly 
\begin{equation}\label{tempo}
(\e_{i+1} \e_i \eta)(t) = (\e_i \eta)(t) + b_{i+1}(t)\an_{i+1}
\end{equation}
etc. Since $\an_k(\can_i) =0 \; \forall k \geq i+2$, we get
\begin{equation}\label{firstone}
\pi_1(t)(\can_i)= \eta(t)(\can_i) + 2 b_i(t) - b_{i+1}(t)
\end{equation}

\noindent
We need to show that the left hand side is $\geq 0$ for all $t$. We
will first show that it is $\geq -1 \; \forall t \in [0,1]$. In fact we claim:
\begin{equation}\label{eqa}
\eta(t)(\can_i) + 2b_i(t) \geq 0 \; \forall t \in [0,1]
\end{equation}
To prove this observe that $\eta(t)(\can_i) + 2b_i(t) =
(\eta(t)+b_i(t)\an_i)(\can_i) = (\e_i \eta)(t)(\can_i)$. By
Proposition~(\ref{eifacts}), (3) we have 
  $\min_t (\e_i \eta)(t)(\can_i) = \min_t \eta(t)(\can_i)   + 1 $. 
Equation~(\ref{eqa}) would thus follow if we show that
\begin{equation}\label{nineteen}
\min_t \eta(t)(\can_i)   = -1
\end{equation}

\noindent
But  $\eta(t)(\can_i) = \pp(t)(\can_i) + \sum_{k=l+1}^{i-2} b_k(t)
\an_k(\can_i) + b_{i-1}(t)(-1)$ where the $b_k$ are increasing functions
with $b_k(0)=0$ and $b_k(1)=1$. Now, $\pp(t)(\can_i) \geq 0$ by 
(\ref{prop-one}), $b_{i-1}(t) \leq 1$ and the intermediate terms in
the sum are 0 since $\an_k(\can_i)$ for $k \leq i-2$. Thus 
$\eta(t)(\can_i) \geq -1$. In fact, $\eta(1)(\can_i) = -1$ since 
$\pp(1)(\can_i) = 0$. This proves Equation~(\ref{nineteen}) and hence
Equation~(\ref{eqa}). Looking back at Equation~(\ref{firstone}), this
means that $\pi_1(t)(\can_i) \geq -1$. Our next step is to show that 
$\pi_1(t)(\can_i)$ never attains the  value -1 for any $t$.

To see this, suppose 
$\pi_1(t_0)(\can_i) = -1$. Then we must have that
\begin{subequations}
\begin{equation}\label{pq}
\eta(t_0)(\can_i) + 2b_i(t_0) =0 
\end{equation}
\begin{equation}\label{rs}
b_{i+1}(t_0) = 1
\end{equation}
\end{subequations}
We look more closely at Equation~(\ref{rs}). By Equation~(\ref{ei}),
$b_{i+1}(t)$ is determined by the values of the function
$(\e_i \eta)(t)(\can_{i+1})$. It is easy to see that if we replace
$\eta$ by $\e_i \eta$ and $\can_i$ by $\can_{i+1}$
 in Equation~(\ref{nineteen}),
then it still holds ({\em the proof is similar}). We record this as
\begin{equation}\label{AB}
\min_t (\e_i \eta)(t)(\can_{i+1})   = -1
\end{equation}
By the definition of $b_{i+1}(t_0)$ (Equation~(\ref{ei}))
 and Equation~(\ref{AB}), there must
exist
$s$, $0 \leq s \leq t_0$ such that $(\e_i\eta)(s) (\can_i) = -1$ i.e,
\begin{equation}\label{XY}
\eta(s)(\can_{i+1}) - b_i(s) = -1
\end{equation}
But observe that the first term on the left hand
side is $\geq 0$. This is because 
$\eta(s)(\can_{i+1}) = \pp(s)(\can_{i+1}) + \sum_{k=l+1}^{i-1} b_k(s)
\an_k(\can_{i+1}) =  \pp(s)(\can_{i+1}) \geq 0$ by Equation~(\ref{prop-one}). 
So the only way Equation~(\ref{XY}) can hold is 
if $\eta(s)(\can_{i+1}) = 0$ and $b_i(s) =1$. Since $b_i$ is an
increasing function and $s \leq t_0$, this means that $b_i(t_0)=1$
as well. Substituting in Equation~(\ref{pq}) we get
$\eta(t_0)(\can_{i}) = -2$. This clearly contradicts
  Equation~(\ref{nineteen}). We have thus shown that
$$\pi_1(t)(\can_i) > -1 \; \forall t \in [0,1]$$
But $\pi_1$ being an L-S path, has the integrality property
(Proposition~(\ref{eifacts}), (5)). Thus $\min_t \pi_1(t)(\can_i) \geq 0$
proving Fact~(\ref{prop-four}). $\hfill{\Box}$

We have thus shown that the path $\pi_1$ is an L-S path, which
satisfies Equations~(\ref{prop-three}) and (\ref{prop-four}). The
situation is now
analogous to the path $\pi$ which satisfies
Equations~(\ref{prop-two}) and (\ref{prop-one}). So we can repeat all
the arguments that came between Equation~(\ref{prop-two}) and
Equation~(\ref{prop-four}) replacing $\pi$ with $\pi_1$ and $U_1$ with
$U_2$ throughout. We thus obtain an L-S path $\pi_2$ of shape
$\LN$ which satisfies

\begin{equation}\label{fin1}
\pi_2(1) - \pi_1(1) \geq \sum_{i=l+2}^{n-r-1} \an_i 
\end{equation}

\begin{equation}
\begin{array}{ll}
\pi_2(t)(\can_i) &= 0 \, \forall i \in U_3 \\
\pi_2(t)(\can_i) &\geq 0 \; \forall t \in [0,1], \, \forall i \in U_3
\end{array}
\end{equation}
where $U_3 :=\{ i: l+2 < i < n-r-1\}$.

It is clear that this process has to stop before the $s^{th}$ stage
 where $s = \dep(\lambda + \mu - \nu)$. 
 This is because the coefficient of $\an_i$ for $(l+s) < i < (n-r+1-s)$ in 
$\LN - \pi_k(1)$ decreases by at least 1 each time $k$ increases by 1 (by 
Equations~(\ref{fin}), (\ref{fin1}), etc). To start with however, we
 know that $\LN - \pi(1) = \gamma\rn$, which is given by
Equation~(\ref{piqis}). Thus the coefficient of these $\an_i$ is $`s'$ to
 begin with. By Proposition~(\ref{eifacts}), (6) the coefficient of
 these $\an_i$ in $\LN - \pi_k(1)$ must be $\geq 0$, forcing 
$k \leq s$. In fact one can show that $k$ must equal $s$,
but we  will not need this fact.

Let $\pi_k$ denote the last path in the list. Then clearly $\LN -
\pi_k(1) = \sum_{i=1}^n c_i \an_i$ with $c_i =0 $ for $i \in V\rn$.
Recall here that
$V\rn =\{l+s < i < n-r+1-s\}$. We
can thus write  $\pln = \e_{i_p} \cdots \e_{i_2} \e_{i_1} \pi_k$ for some
$i_1,i_2,\cdots,i_p \not\in V\rn$. In summary , if we define $\pi_0 =
\pi$ and $\pi_{k+1} = \pln$, then for each $0 \leq j \leq k$, we
have shown that $\pi_j$ can be obtained from $\pi_{j+1}$ by repeated action
of elements of $T = \{ f_i\rn: i \not\in V\rn\} \cup \{g\rn,h\rn\}$. Thus $\pi_0 =
\pi$ can be obtained from $\pi_{k+1} = \pln$ by the action of elements
of $T$. Hence $\pi \in \mathcal{P}^0(\lambda,\mu,\nu,n)$.
This concludes the proof of Proposition~(\ref{thecrux}).$\hfill{\Box}$

\subsection{Special cases}

We now restrict $\lambda, \mu, \nu$ to be certain special
 types of double-headed
 weights and say something more about the stable multiplicities for
 these types.

\begin{enumerate}
\item
{The first type we shall consider was the starting point
for this present work. Let $\lambda$, $\mu$ and $\nu$ be {\em
  single headed} dominant weights supported on the ``tail'' portion
i.e, $\lambda = (0,x)$, $\mu = (0,y)$, $\nu = (0,z)$ for $x,y,z \in
\mathcal{H}_1^+$. If $x=(x_1,x_2,x_3,\cdots)$, $y=(y_1,y_2,y_3,\cdots)$
$z=(z_1,z_2,z_3,\cdots)$, then $|\lambda|_X  + |\mu|_X = |\nu|_X$
implies that $\sum_i i\,(x_i+y_i) = \sum_i i \,z_i$. So, when thought
of as partitions, the numbers of boxes in the Young diagrams of $x$
and $y$ add up to that of $z$. It is clear from the proof of
Proposition~(\ref{piqis}) (or alternatively from
Equation~(\ref{tailmod})) that in this case, there exist integers 
$q_1, q_2, \cdots, q_{r-1}$ such that 
\begin{equation}\label{lmn}
\LN + \mn - \nn = \sum_{i=1}^{r-1} q_i \an_{n-i+1} \;\; \forall n \geq r
\end{equation}
(the $p_i$ and $s$ in Proposition~(\ref{piqis}) are 0). Here $r =
\max(\supp(x),\supp(y),\supp(z))$. 

Let us now consider the Dynkin diagrams $A_n$. Clearly $x,y$ and $z$
define dominant weights of $A_n$ by setting
$$x\rn := \sum_{i=1}^{\supp(x)} x_i \overline{\omega}_i\rn$$
where $\overline{\omega}_i\rn$ denotes the fundamental weight of $A_n$
corresponding to the node $n-i+1$. The definitions of $y\rn$ and
$z\rn$ are similar. It is also clear that
\begin{equation}\label{xyz}
x\rn + y\rn -z\rn = \sum_{i=1}^{r-1} q_i \tilde{\alpha}\rn_{n-i+1} \;\; \forall n \geq r
\end{equation}
where $r$ and $q_i$ are the same integers as above, and the
$\tilde{\alpha}\rn_{i}$ denote the simple roots of $A_n$ (as opposed
to $X_n$).

By Littelmann's tensor product decomposition  formula, $\clmn(n)$ is
the number of $\mn$ dominant L-S paths $\eta = f_{i_k}\rn
f_{i_{k-1}}\rn \cdots f_{i_1}\rn(\pln)$ which satisfy 
$\sum_{j=1}^k \alpha_{i_j}\rn = \sum_{i=1}^{r-1} q_i
\an_{n-i+1}$. Similarly, if
  $\tilde{c}_{xy}^{\, z}(n) :=$ multiplicity of $L(z\rn)$ in
 $L(x\rn) \otimes L(y\rn)$ as representations of $\mathfrak{g}(A_n)$, then
$\tilde{c}_{xy}^{\, z}(n)$ is the number of L-S paths
$\tilde{\eta}$  (in $\mathfrak{h}^*(A_n)$) of the form
 $\tilde{f}_{i_k}\rn \tilde{f}_{i_{k-1}}\rn \cdots 
\tilde{f}_{i_1}\rn(\pi_{x\rn})$ which are $y\rn$ dominant and satisfy 
$\sum_{j=1}^k \tilde{\alpha}_{i_j}\rn = \sum_{i=1}^{r-1} q_i 
\tilde{\alpha}\rn_{n-i+1}$. The $\tilde{f}_i\rn$ denote lowering
operators of $A_n$.

\vspace{0.1in}
Observe that the lowering operators involved in both cases correspond
to the rightmost $r-1$ nodes. The above expressions for tensor product
multiplicities for $X_n$ and $A_n$ clearly imply that $\clmn(n) =
\tilde{c}_{xy}^{\, z}(n)$. Taking $n$ large enough,
Theorem~(\ref{mainthmrestated}) implies that $\clmn(n) =
\clmn(\infty)$. By the classical theory for $A_n$, 
we know that $\tilde{c}_{xy}^{\, z}(n)$ equals the
Littlewood-Richardson coefficient $LR_{x,y}^{\,z}$ corresponding to
$x,y,z$ (considered as partitions). Thus $\clmn(\infty) = LR_{x,y}^{\,z}$.

In summary, as long as all three dominant weights of $X_n$ under consideration
are supported on the ``tail'', their behavior is exactly like 
dominant weights of $A_n$.}

\item { Let us consider a different situation. Let $\lambda, \mu$ be single
  headed weights, supported on the tail  as above, but now let $\nu$ be a
  single headed weight supported on the ``head'' i.e, near the $X$
  portion of the Dynkin diagram. So $\lambda = (\zero,x)$, $\mu = (\zero
  , y)$, $\nu = (z,\zero)$ for $x,y,z \in \mathcal{H}_1^+$.

\vspace{0.05in}
If $X_n = A_n$, then $\clmn(\infty) = LR_{\lambda,\mu}^{\,\nu} =0$,
since the Littlewood-Richardson rule implies that any $\nu$ for which 
 $LR_{\lambda,\mu}^{\,\nu} \neq 0$ has to be supported on the rightmost
$k$ nodes, where $k = \supp(x) + \supp(y)$. For large $n$, $\nu =
(z,\zero)$ fails to meet this criterion. One can also obtain this fact
from our point of view. Notice that $|\nu|_A = |(z,\zero)|_A=
 \sum_{i=1}^{\supp(z)} i\,z_i > 0$, while $|\lambda|_A + |\mu|_A =
 |(\zero,x+y)|_A = - \sum_i i\, (x_i + y_i) < 0$. Hence $|\nu|_A \neq 
|\lambda|_A + |\mu|_A$. Proposition~\eqref{numofboxes} implies 
$\clmn(\infty) = 0$.

\vspace{0.1in}
However for $X_n = E_n$, $\clmn(\infty)$ could be positive, as we saw
in Example~\eqref{eg2} for
$x=y=z=\epsilon_1$. Observe that the contradiction obtained above for
$A_n$ in terms of the number of boxes function disappears for
$E_n$. From the  definition, it follows that $|\lambda|_E + |\mu|_E >
0$ while $|\nu|_E = \sum_i a_i z_i = 2z_1 + 4z_2 + \cdots + (10-k)
z_k$ which is positive for many choices of $z_i$.

In summary, for $A_n$, 
if $\lambda$ and $\mu$ are single headed and supported on
the tail portion, any $\nu$ for which $\clmn(\infty) > 0$ must also be
single headed and supported on the tail portion. However for $E_n$,
this is not the case. In fact, there can even exist a $\nu$, 
supported on the ``head'' portion for which $\clmn(\infty) > 0$. In
a sense, information that is localized at one end (the tail) of the Dynkin
diagram of $E_n$ propagates to the other end.
}
\end{enumerate}

\section{The Stable Representation Ring}

Having established that the multiplicities $\clmn(n)$ stabilize, we
shall now use the stable values $\clmn(\infty)$ as structure constants
to define a multiplication operation $*$ on a space $\srr$. We
shall call $\srr$ the {\em stable representation ring of type
  $X$}. 

In type $A$, 
the associativity of $*$ will follow directly from
the associativity of the tensor product. But for general type $X$,
using $\clmn(\infty)$ as structure constants means that we only keep
the stable terms in the tensor product decomposition and discard the
``transient'' ones. Associativity of $*$ is no longer obvious. The
goal of this section is to show that associativity still holds and
that $\srr$ becomes a genuine $\complex$ - algebra.

We assume $X$ is an extensible marked Dynkin diagram with $d$ nodes.
We shall consider
the tensor product of three or more irreducible integrable highest weight 
representations and study its decomposition. First we will need a
technical lemma concerning the large $n$ behavior of the set of 
dominant weights $P^+(X_n)$. We prove this so called Interval
Stabilization lemma in Section~(\ref{issection}) and then use it in
Section~(\ref{appl}) to look at stable multiplicities in $k$-fold
tensor products . The stable representation ring will be defined in
 Section~(\ref{stabrepring}).

\subsection{Interval Stabilization}\label{issection}
First, some notation that will be needed to state our lemma:
Let $\lambda_1 = (x,y)$ , $\lambda_2=(z,w) \in \mathcal{H}_2$ be such that 
$|\lambda_1|_X = |\lambda_2|_X$. Let $l = \max(d,\ell(x),\ell(z))$
and $r = \max(\ell(y),\ell(w))$. Proposition~(\ref{piqis}) implies that
there exist integers $p_i \; (1\leq i \leq l-1$) , $q_j \; (1\leq j \leq r-1)$
and $s$ such that for $n \geq l+r$
\begin{equation}\label{piqiagain}
\lambda_1^{(n)}-\lambda_2^{(n)} = \sum_{i=1}^{l-1} p_i \alpha_i^{(n)} +
\sum_{i=l}^{n-r+1} s \alpha_i^{(n)} + 
\sum_{i=n-r+2}^{n} q_{(n-i+1)} \alpha_i^{(n)}
\end{equation}
We define a partial order $\geq$ on $\mathcal{H}_2$ by requiring that 
$\lambda_1 \geq \lambda_2$ iff $|\lambda_1|_X = |\lambda_2|_X$ and the
$p_i, q_j, s$ which occur in Equation~(\ref{piqiagain}) are all
non-negative.

It is easy to check that $\geq$ is a partial order on $\mathcal{H}_2$
and that $\lambda_1 \geq \lambda_2$ implies that 
$\lambda_1+\mu \geq \lambda_2+\mu$. We also have these equivalent
conditions which follow from the arguments of Section \ref{nob}:
\begin{equation*}
\begin{array}{ll}
\lambda_1 \geq \lambda_2 &\Leftrightarrow \lambda_1\rn - \lambda_2\rn
\in Q^+(X_n) \; \forall \mbox{ large } n \\
&\Leftrightarrow \lambda_1\rn - \lambda_2\rn
\in Q^+(X_n) \; \mbox{ for infinitely many values of  } n\\
&\Leftrightarrow |\lambda_1|_X = |\lambda_2|_X \text{ and } 
\lambda_1\rn - \lambda_2\rn
\in Q^+(X_n) \; \text{ for some value of  } n \geq l+r
\end{array}
\end{equation*}
Recall that the usual partial order $\geq$ on $\mathfrak{h}^*(X_n)$ is
defined by $\beta \geq \beta^{\prime}$ iff $\beta - \beta^{\prime} \in
Q^+(X_n)$ ($\beta, \beta^{\prime} \in \mathfrak{h}^*(X_n)$). Hence 
for $\lambda_1, \lambda_2 \in \mathcal{H}_2$, $\lambda_1 \geq
\lambda_2$ iff $\lambda_1\rn \geq \lambda_2\rn$ in
$\mathfrak{h}^*(X_n)$ for all large $n$.
We now state our main lemma.
\begin{lemma}\label{is}
(Interval Stabilization) Let  $\lambda_1, \lambda_2 \in
  \mathcal{H}_2^+$ with $\lambda_1 \geq \lambda_2$. Let $\il := \{ \gamma
  \in \mathcal{H}_2^+ : \lambda_1 \geq \gamma \geq \lambda_2\}$ and 
$\INL := \{ \beta  \in P^+(X_n) : \lambda_1\rn \geq \beta \geq
  \lambda_2\rn\}$ for $n$ larger than the lengths of $\lambda_1$ and
  $\lambda_2$. Then 
\begin{enumerate}
\item $\il$ is a finite set
\item There exists $N$ such that for all $n \geq N$, $\INL = \{
  \gamma\rn: \gamma \in \il \}$
\end{enumerate}
\end{lemma}

\noindent
{\bf Proof:} Let $\lambda_1 = (x,y)$ , $\lambda_2=(z,w)$. 
Let $l = \max(d,\ell(x),\ell(z))$ and $r =
\max(\ell(y),\ell(w))$. We had by Equation~(\ref{piqiagain})
\begin{equation*}
\lambda_1^{(n)}-\lambda_2^{(n)} = \sum_{i=1}^{l-1} p_i \alpha_i^{(n)} +
\sum_{i=l}^{n-r+1} s \alpha_i^{(n)} + 
\sum_{i=n-r+2}^{n} q_{(n-i+1)} \alpha_i^{(n)}
\end{equation*}
Since $\lambda_1 \geq \lambda_2$, the $p_i$, $q_j$ and $s$ are all
non-negative.  
Set $N = l+r + 2s$. Fix $n \geq N$. Let $\overline{U}\rn := \{ l \leq
i \leq n-r+1\}$, $U\rn :=  \{ l < i < n-r+1\}$, $\overline{V}\rn :=
\{l+s \leq i \leq n-r+1-s\}$ and $V\rn := \{ l+s < i < n-r+1-s\}$.

Pick $\beta \in \INL$ i.e, $\beta \in P^+(X_n)$ and $\LN_1 \geq \beta
\geq \LN_2$. Hence
\begin{equation}\label{inequality}
0 \leq \beta - \LN_2 \leq \LN_1 - \LN_2
\end{equation}
If $ \beta - \LN_2 = \sum_{i=1}^n b_i \alpha_i\rn$,
Equations~(\ref{inequality}) and (\ref{piqiagain}) imply that 
$$(i) \;\; 0 \leq b_i \leq s \; \forall i \in \ou$$
 Since $\beta \in P^+(X_n)$ and
$\LN_2(\can_i) =0 \; \forall i \in U\rn$, we have $(\beta -
\LN_2)(\can_i) \geq 0 \; \forall i \in U\rn$. This gives us the
following additional condition on the $b_i 's$
$$ (ii) \;\; 2b_i - b_{i-1} - b_{i+1} \geq 0 \; \forall i \in U\rn$$

\noindent
{\bf Claim:} $b_i$ is a constant on $\ov$ i.e, $b_i = b_j \; \forall
i,j \in \ov$. {\em Proof:} Suppose not, then there exists $i$ such
that  $i, i+1 \in \ov$, but $b_i \neq b_{i+1}$. {\bf Case 1:}
 Suppose $b_i >
b_{i+1}$. Condition (ii) implies that $b_{i+2} \leq 2b_{i+1} - b_i <
b_{i+1}$. Similarly we conclude $b_{i+3} < b_{i+2}$ etc. So we have a
strictly descending sequence $b_i > b_{i+1} > b_{i+2} > \cdots >
b_{n-r+1}$. The number of terms in this sequence is 
$(n-r-i+2) \geq s+2$ (since
$i+1 \in \ov$ means that $i+1 \leq n-r+1-s$) and by (i) we know that
each term in the sequence lies between $0$ and $s$ . This is a clear 
contradiction. {\bf Case 2:} Suppose $b_i < b_{i+1}$. We proceed as above to
conclude that $b_l < b_{l+1} < \cdots < b_i < b_{i+1}$. The number of
terms in this  ascending sequence is $(i-l+2) \geq s+2$ 
(since $i \in \ov$ implies $i \geq
l+s$). Again a contradiction.$\hfill{\Box}$

We denote the constant value by $k$. Hence $k = b_i \; \forall i \in
\ov$.

\noindent
{\bf Consequences:} (1) $(\beta - \LN_2)(\can_i) = 0 \; \forall i \in
V\rn$. This is clear since the left hand side is just $2b_i -
b_{i-1}-b_{i+1}$ and $i,i-1,i+1 \in \ov$. 

Since $\LN_2(\can_i) = 0 \; \forall i \in V\rn$, this also means that 
$\beta(\can_i) = 0 \; \forall i \in V\rn$. This implies that if $\gamma
= (t,u) \in \il$, then $\max(d,\ell(t)) \leq l+s$ and $\ell(u) \leq
r+s$, since $\gamma^{(m)} \in I^{(m)}(\lambda_1,\lambda_2)$
 for all large $m$. 
Since $n \geq l+r+2s$, we get a well defined, 
injective map $\phi_n : \il \rightarrow \INL$ defined by
$\phi_n(\gamma) := \gamma\rn$. Since $\INL = \{ \beta \in P^+(X_n):
\LN_1 \geq \beta \geq \LN_2\}$ is a finite set, $\il$ must be finite
too. This proves statement (1) of the Lemma.

(2) Since $\beta(\can_i) = 0 \; \forall i \in V\rn$, we can write 
$$\beta = \sum_{i=1}^{l+s} c_i \omega_i\rn + \sum_{j=1}^{r+s} 
\tilde{c}_j \overline{\omega}_j\rn$$ 
Define $c := (c_1, c_2, \cdots, c_{l+s}, 0,0,\cdots)$, 
$\tilde{c} := (\tilde{c}_1, \tilde{c}_2, \cdots, \tilde{c}_{r+s},
0,0,\cdots)$ and set $\gamma:=(c, \tilde{c}) \in
\mathcal{H}_2^+$. Then $\beta = \gamma\rn$. If we show that
$|\gamma|_X = |\lambda_2|_X = |\lambda_1|_X$, then $\gamma$ would be
an element of $\il$ since $\gamma\rn \in \INL$. This would prove statement
(2) of the Lemma as well.

We know that 
\begin{equation}\label{k}
\gamma\rn - \LN_2 = \sum_{i=1}^{l+s-1} b_i \an_i + k
\sum_{i=l+s}^{n-r-s+1} \an_i + \sum_{j=1}^{r+s-1} \tilde{b}_j\an_{n-j+1}
\end{equation}
where for $1 \leq j \leq r+s-1$, $\tilde{b}_j := b_{n-j+1}$. Since
$\gamma$ is supported on the first $l+s$ and last $r+s$ nodes, it is
clear from Equation~(\ref{k}) above that for all $m \geq n$,
$\gamma^{(m)} - \lambda_2^{(m)}$ is given by 
\begin{equation*}
\gamma^{(m)} - \lambda_2^{(m)} = \sum_{i=1}^{l+s-1} b_i \alpha_i^{(m)} + k
\sum_{i=l+s}^{m-r-s+1} \alpha_i^{(m)} + 
\sum_{j=1}^{r+s-1} \tilde{b}_j\alpha_{m-j+1}^{(m)}
\end{equation*}
obtained by ``elongating'' the string of $k$'s in the
middle. Hence $\gamma^{(m)} \equiv \lambda_2^{(m)} \pmod{Q(X_m)}$
for all $m \geq n$. By the arguments of Section \ref{nob},
this implies that $|\gamma|_X = |\lambda_2|_X$. This finishes the
proof of Lemma~(\ref{is}) $\hfill{\Box}$

\subsection{$k$-fold tensor products}\label{appl}

To extend our main theorem~\eqref{mainthm}, we now turn to tensor
products of three or more irreducible representations. We ask if
multiplicities in $k$-fold tensor
products also stabilize. We shall first show that this remains true. 
Secondly, it is not
obvious that one can understand stable multiplicities in $k$-fold
tensor products by understanding stable multiplicities in successive
binary tensor products. Happily it turns out that this can also be done.

\begin{definition}
Let $\lambda_1,\ldots,\lambda_k$ and $\nu\in\mathcal{H}_2^+$
 be double-headed
weights.  Define $\clarge (n)$ to be
the multiplicity of the representation $L(\nu^{(n)})$ in the $k$-fold
tensor product $L(\lambda_1^{(n)})\otimes\cdots\otimes L(\lambda_k^{(n)})$.
If this is independent of $n$ when $n$ is large, let
$\clarge (\infty)$ denote its stable value.
\end{definition}

This generalizes the preceding use of $\clmn(n)$.

\begin{theorem}\label{A}
If $|\lambda_1|_X+\cdots+|\lambda_k|_X=|\nu|_X$, then
$\clarge (n)$ is indeed independent of $n$
for $n$ sufficiently large.  Moreover, the stable value is related to
the stable multiplicities in successive binary tensor products in the
usual way:
\begin{equation}\label{kvstwo}
\clarge(\infty) =
\sum_{\mu_1,\ldots,\mu_{k-2}\in\mathcal{H}_2^+} 
 c_{\lambda_1\lambda_2}^{\,\mu_1}(\infty)
 c_{\mu_1\lambda_3}^{\,\mu_2}(\infty)
 \cdots
 c_{\mu_{k-3}\lambda_{k-1}}^{\,\mu_{k-2}}(\infty)
 c_{\mu_{k-2}\lambda_k}^{\,\nu}(\infty)
\end{equation}
\end{theorem}
If $n$ is finite, then equation~\eqref{kvstwo} is clearly true, if we
replace the $\infty$'s by $n$ and let the sum range over all $\mu_i
\in P^+(X_n)$. This holds since
$$L(\lambda_1\rn) \otimes \cdots \otimes L(\lambda_k\rn) =
(\ldots((L(\lambda_1\rn) \otimes L(\lambda_2\rn)\,)\otimes
L(\lambda_3\rn)\,)\otimes \ldots \otimes  L(\lambda_k\rn)\,)$$
We will use the Interval stabilization lemma~\eqref{is} to show that
when $n$ is large enough, then the ranges of $\mu_i$ we must sum over
also stabilize. This will prove both parts of Theorem~\eqref{A}.

\noindent
{\bf Proof of Theorem~\eqref{A}:} The essence of the proof is the case
$k=3$. The general case follows by making modifications in the obvious
places.

Consider $\csmall(n)$.  Since 
\begin{equation}\label{tpassoc}
 L(\lambda_1\rn) \otimes  L(\lambda_2\rn)  \otimes  L(\lambda_3\rn)
\cong 
 (L(\lambda_1\rn) \otimes  L(\lambda_2\rn)\,)  \otimes  L(\lambda_3\rn)
\end{equation}
we have
\begin{equation}\label{tripdoub}
\csmall(n)  = \sum_{\beta \in P^+(X_n)}
c_{\lambda_1\rn,\,\lambda_2\rn}^{\,\beta} \,\cdot
c_{\beta,\,\lambda_3\rn}^{\,\nu\rn}
\end{equation}
By a  mild abuse of
 notation, we let $c_{\beta_1,  \beta_2}^{\,\beta_3}$ 
denote the multiplicity of $L(\beta_3)$ in
  $L(\beta_1) \otimes L(\beta_2)$ (all representations of $X_n$),
for $\beta_i \in P^+(X_n),\, i=1,2,3$. 
Now  if
\begin{equation}\label{ineq}
 c_{\lambda_1\rn,\,\lambda_2\rn}^{\,\beta} > 0 \mbox{ and } 
c_{\beta,\,\lambda_3\rn}^{\,\nu\rn} > 0
\end{equation}
we get $\lambda_1\rn + \lambda_2\rn \geq \beta$ and  
$\beta + \lambda_3\rn \geq \nu\rn$. Hence
$$ \lambda_1\rn + \lambda_2\rn + \lambda_3\rn \geq \beta +
 \lambda_3\rn \geq \nu\rn$$
We note that  $\sum_{i=1}^3 \lambda_i\rn \geq \nu\rn$
 together with  $\sum_{i=1}^3 |\lambda_i|_X =|\nu|_X$ implies
 that $\sum_{i=1}^3 \lambda_i \geq \nu$ in the partial
 order on $\mathcal{H}_2^+$.

Let $\tilde{\beta} = \beta + \lambda_3\rn$.
We can  now apply the Interval Stabilization
Lemma~(\ref{is}). This gives us an integer $N^{\prime}$ such that for 
$n \geq N^{\prime}$, $\tilde{\beta} = \gamma\rn$ for some $\gamma \in
 I(\lambda_1+\lambda_2+\lambda_3, \nu)$. 
So $\beta = \gamma\rn - \lambda_3\rn$. Let 
$$F_{\lambda_3} 
 := \{ \gamma - \lambda_3: \gamma \in I(\lambda_1+\lambda_2+\lambda_3,
 \nu)\} \cap \mathcal{H}_2^+$$
The only possible solutions $\beta$ to (\ref{ineq}) are
 $\beta = \delta\rn$, for 
$\delta \in F_{\lambda_3}$.
Thus
$$\csmall(n) = \sum_{\delta \in 
F_{\lambda_3}}
  c_{\lambda_1\lambda_2}^{\,\delta}(n)\,c_{\delta\lambda_3}^{\nu}(n)$$
Since the number of terms in this sum is finite, we can pick $N \geq
  N^{\prime}$ such that for all $n \geq N$ and all 
$\delta \in F_{\lambda_3}$,
$c_{\lambda_1\lambda_2}^{\,\delta}(n) = 
c_{\lambda_1\lambda_2}^{\,\delta}(\infty)$ 
and $c_{\delta\lambda_3}^{\,\nu}(n) = 
c_{\delta\lambda_3}^{\,\nu}(\infty)$. Hence for
  all $n,m \geq N$, $\csmall(n) = \csmall(m)$. 
We've thus shown that
  the multiplicities of representations in the triple tensor product
 do stabilize. We've in fact also shown:
\begin{equation}\label{tri}
\csmall(\infty) =\sum_{\delta \in F_{\lambda_3}}
  c_{\lambda_1\lambda_2}^{\,\delta}(\infty) \, 
c_{\delta\lambda_3}^{\nu}(\infty) \,
=\sum_{\gamma \in \mathcal{H}_2^+}  c_{\lambda_1\lambda_2}^{\,\gamma}(\infty)
\, c_{\gamma\lambda_3}^{\nu}(\infty)
\end{equation}
For the last equality, observe by usual arguments that 
$c_{\lambda_1\lambda_2}^{\,\gamma}(\infty) > 0$ and 
$c_{\gamma\lambda_3}^{\,\nu}(\infty) > 0$ imply that $\lambda_1 +
\lambda_2 \geq \gamma$ and $\gamma + \lambda_3 \geq \nu$ 
in the partial order on
$\mathcal{H}_2$. Hence $\gamma + \lambda_3 \in I(\lambda_1 + \lambda_2
+ \lambda_3, \nu)$. 
So $\gamma \in \{\delta - \lambda_3 : \delta \in I(\lambda_1 +
\lambda_2 + \lambda_3, \nu)\} \cap  \mathcal{H}_2^+ = 
F_{\lambda_3}$. $\hfill{\Box}$

\begin{remark}
In the above proof, instead of \eqref{tpassoc} we could have started
from the fact that 
 $L(\lambda_1\rn) \otimes  L(\lambda_2\rn)  \otimes  L(\lambda_3\rn)
\cong 
 L(\lambda_1\rn) \otimes  (L(\lambda_2\rn)  \otimes
 L(\lambda_3\rn))$. It is clear that we would have obtained the
 following equation analogous to equation~\eqref{tri}:
\begin{equation}\label{retri}
\csmall(\infty) =\sum_{\delta \in F_{\lambda_1}}
  c_{\lambda_1\delta}^{\,\nu}(\infty) \, 
c_{\lambda_2\lambda_3}^{\delta}(\infty)
=\sum_{\gamma \in \mathcal{H}_2^+}  c_{\lambda_1\gamma}^{\,\nu}(\infty)
\, c_{\lambda_2\lambda_3}^{\gamma}(\infty)
\end{equation}
\end{remark}
\subsection{The stable representation ring $\srr$}\label{stabrepring}

Theorem~\eqref{A} is  key to our definition of
$\srr$. First let  $\mathcal{R}$ denote the $\complex$
vector space with basis $\{v_{\lambda}: \lambda \in \mathcal{H}_2^+\}$
and $\widehat{\mathcal{R}}$ be its formal completion i.e, 
 $\widehat{\mathcal{R}}$ is the set
$$ \left\{ \sum_{\lambda \in \mathcal{H}_2^+} c_{\lambda} v_{\lambda} \,:\,
  c_{\lambda} \in \complex \right\}$$
 of all formal infinite series in the $v_{\lambda}$.

We define a multiplication operation on the basis elements $v_{\lambda}$.
$$v_{\lambda} * v_{\mu} := \sum_{\gamma \in \mathcal{H}_2^+}
c_{\lambda\mu}^{\,\gamma}(\infty) \,v_{\gamma}$$
Equation~(\ref{tri}) shows that $(v_{\lambda} * v_{\mu})*v_{\nu}
:= \sum_{\pi \in  \mathcal{H}_2^+} \left( \sum_{\gamma \in
\mathcal{H}_2^+}   c_{\lambda\mu}^{\,\gamma}(\infty)\,
 c_{\gamma\nu}^{\pi}(\infty) \right) v_{\pi} $ 
 is equal to $\sum_{\pi} \cp(\infty) v_{\pi}$ and hence well defined.
Analogously, equation~\eqref{retri} guarantees that
  $v_{\lambda} * (v_{\mu} * v_{\nu})$ is also well defined and equal to 
  $\sum_{\pi} \cp(\infty) v_{\pi}$. 
Thus:
\begin{equation}\label{assoc}
(v_{\lambda} * v_{\mu})*v_{\nu} = v_{\lambda} * (v_{\mu} * v_{\nu})
\end{equation}
Looking back on section~\eqref{appl}, we see that this associativity
is essentially a consequence of the associativity of the tensor
product:
$$(L(\lambda) \otimes L(\mu)) \otimes L(\nu)
  \cong L(\lambda) \otimes (L(\mu) \otimes L(\nu))$$
Further,  theorem~\eqref{A} on $k$-fold tensor products shows that the
  product $v_{\lambda_1} * v_{\lambda_2} * \cdots * v_{\lambda_k}$ of
  finitely many $v_{\lambda_i}$'s is necessarily well defined, since
  it is equal to 
$\sum_{\nu \in \mathcal{H}_2^+} \clarge(\infty) v_{\nu}$. We then make
  the following definition:
\begin{definition}
Let $\srr$ denote the subspace of
$\widehat{\mathcal{R}}$ spanned by the set 
$$\{ v_{\lambda_1} * v_{\lambda_2} * \cdots * v_{\lambda_k}: k \geq 0,
\lambda_i \in \mathcal{H}_2^+\}$$
consisting of all finite products of the $v_{\lambda}$'s. 
\end{definition}
Clearly $\srr$ is an associative, commutative $\complex$ algebra
with respect to the operations of adddition and $*$. We call
$\srr$ the {\em stable representation ring of type $X$}.
One thinks of $\srr$ as encoding information about how tensor
products decompose as $n \rightarrow \infty$. 

When $X$ is of type $A$, $\srr[A]$ can be identified with the
polynomial algebra $\complex[x_1,y_1,x_2,y_2,\cdots]$ via the map that
sends $x_i \mapsto v_{(\epsilon_i,0)}$ and $y_i \mapsto  v_{(0,\epsilon_i)}$.
Here $\epsilon_i$ denotes the element $(0,0,\cdots,1,0,\cdots) \in
\mathcal{H}_1^+$ with the 1 in the $i^{th}$ place. If we introduce
$\integers$-gradations on these two algebras,
by setting $\deg(x_i) = i = -\deg(y_i)$ and $\deg(v_{\lambda}) =
|\lambda|_A$ for $\lambda \in \mathcal{H}_2^+$, then the above map
defines an isomorphism of {\em graded} algebras.

Equivalently, one can view $\srr[A]$ as the tensor product of two
copies of the ring of symmetric functions by identifying $x_i$ and
$y_i$ with the $i^{th}$ elementary symmetric polynomials in the
variables $z_i$ and $w_i$ respectively. Here the gradation would be
$\deg(z_i) = 1 = -\deg(w_i)$.  
In this picture,  the subalgebra of $\srr[A]$ generated by the elements
$\{v_{(\lambda,0)}\,: \lambda \in \mathcal{H}_1^+\}$ is 
isomorphic to the algebra of symmetric functions. Our map above
 sends $v_{(\lambda,0)}$ to the Schur function
 $s_{\lambda}(z_1,z_2,\cdots)$. 

For general $X$, a better understanding
of the structure of $\srr$ might shed more light on the representation
theory of the $X_n$. 
We conclude by mentioning an important open
problem: How far does the ring $\srr$ characterize the series $X_n$ ?
Can there exist an isomorphism $\srr \cong \srr[Y]$ for two different
``types'' $X$ and $Y$ ?

%

\vspace{0.1in}
\noindent
\scriptsize{DEPARTMENT OF MATHEMATICS, BRANDEIS UNIVERSITY, WALTHAM, MA 02454}\\
{\em E-mail address} : {\tt kleber@brandeis.edu}

\vspace{0.1in}
\noindent
{\scriptsize DEPARTMENT OF MATHEMATICS, UNIVERSITY OF CALIFORNIA, BERKELEY, CA 94720}\\
{\em E-mail address} : {\tt svis@math.berkeley.edu}
\end{document}